\documentclass[final,leqno]{siamltex}
\usepackage{graphicx}
\usepackage{amsmath}
\usepackage{amsfonts}
\usepackage{amssymb}
\usepackage{mathrsfs}
\usepackage{verbatim}
\usepackage{subfigure}
\usepackage[notref, notcite, color]{showkeys}
\usepackage{datetime}
\title{Numerical Analysis for a System Coupling 
Curve Evolution to Reaction-Diffusion
on the Curve}%\footnotemark[1]}
\author{John W. Barrett\footnotemark[2] \and 
        Klaus Deckelnick  \footnotemark[3]\ \and 
        Vanessa Styles \footnotemark[4]}%

\newcommand{\bR}{\mathbb{R}}
\newcommand{\matpartv}{\partial_t^\bullet}

\newcommand{\vE}{\vec{E}}
\newcommand{\vX}{\vec{X}}
\newcommand{\vx}{\vec{x}}
\newcommand{\vxi}{\vec{\xi}}

\newcommand{\vtau}{\vec{\tau}}
\newcommand{\vnu}{\vec{\nu}}
\newcommand{\vTau}{\vec{\mathcal{T}}}
\newcommand{\vNu}{\vec{\mathcal{V}}}
\newcommand{\Vh}{V^h}

\newcommand{\Ip}{\mathbb{I}}

\newcommand{\drho}{\; {\rm d}\rho}
\newcommand{\dt}{\; {\rm d}t}

\def\epsilon{\varepsilon}

\begin{document}

\maketitle

\renewcommand{\thefootnote}{\fnsymbol{footnote}}

%\footnotetext[1]{The authors certify that the general content of this manuscript is not submitted, accepted, or published elsewhere, including conference proceedings.}

\footnotetext[2]{Department of Mathematics, 
Imperial College London, London, SW7 2AZ, UK}
\footnotetext[3] {Institut f{\"u}r Analysis und Numerik,
 Otto-von-Guericke-Universit{\"a}t, 39106 
Magdeburg, Germany}
\footnotetext[4]{Department of Mathematics, University of Sussex, Brighton, BN1 9RF, UK}

\begin{abstract}
We consider a finite element approximation for a system
consisting of the evolution of a closed planar curve by forced curve shortening flow
coupled to a reaction-diffusion equation on the evolving curve. 
The scheme for the curve evolution is based on a parametric description 
allowing for tangential motion,
whereas the discretisation  for the PDE on the curve uses an idea from  
\cite{DziukE07}. We prove optimal error bounds for the resulting  fully
discrete approximation and present numerical experiments. These confirm our estimates
and also illustrate the advantage of the tangential motion of the mesh points in practice. 
\end{abstract}



%\pagestyle{myheadings}
%\thispagestyle{plain}
%\markboth{K. DECKELNICK, C.~M. ELLIOTT AND T. RANNER}{UNFITTED FINITE ELEMENT METHODS FOR SURFACE PDES}

\section{Introduction} \label{sec:1}

The aim of this paper is to analyse a fully discrete numerical scheme for approximating a solution
of the following system: find a family of planar, closed curves $(\Gamma(t))_{t \in [0,T]}$  and a function
$w: \bigcup_{t \in [0,T]} \bigl( \Gamma(t) \times \lbrace t \rbrace \bigr) \rightarrow \bR$ such that
\begin{subequations}
\begin{alignat}{2}
v & = \kappa + f(w)   \qquad &&\mbox{on }\Gamma(t), \quad t \in (0,T], \label{Gamma1}\\
\matpartv w  & = d \,w_{ss} + \kappa\,v\,w + g(v,w)  \qquad &&
\mbox{on }\Gamma(t), \quad t \in (0,T], \label{Gamma2}
\end{alignat}
\end{subequations}
subject to the initial conditions
\begin{equation}  \label{Gamma3}
\Gamma(0)=\Gamma^0, \qquad w(\cdot,0) = w^0  \mbox{ on } \Gamma^0.
\end{equation}
Here, $v$ and $\kappa$ are the normal velocity and the curvature of $\Gamma(t)$ corresponding to
the choice $\vnu$ of a unit normal, while $s$ is the arclength parameter on $\Gamma(t)$. Furthermore,
$\matpartv w$ denotes the material derivative of $w$, i.e. $\matpartv w = w_t + v\, 
\frac{\partial w}{\partial \nu}$.  
%$\bar{w}$ is an extension of $w$ into a space-time neighbourhood. 
Finally, $d\in \bR_{>0}$, 
$f:\bR \rightarrow \bR$, $g:\bR \times \bR \rightarrow
\bR$, the closed curve $\Gamma^0$ and $w^0:\Gamma^0 \rightarrow \bR$ are all given. 
The  system (\ref{Gamma1},b) %, (\ref{Gamma2})
couples the evolution of the curves $(\Gamma(t))_{t \in [0,T]}$ 
by forced curve shortening flow to a parabolic PDE on the moving
curves. It occurs for example as a sharp--interface model for diffusion induced grain
boundary motion: in this setting $\Gamma(t)$ represents a grain boundary separating the
crystals of a thin polycrystalline film of metal that is placed in a vapour containing
another metal; atoms from the vapour diffuse into the 
film along the grain boundaries causing them to move. A thorough description of the 
physical set-up can be found in 
 \cite{Hand88}, while an existence and uniqueness result has been obtained  
in \cite{MayerS99}. 

In what follows we shall describe the evolving curves $\Gamma(t)$ 
with the help of a parametrisation $\vx(\cdot,t) : \Ip \rightarrow \bR^2$,
where $\Ip := \bR \setminus \mathbb{Z}$ is the periodic unit interval. Then
\begin{equation}  \label{tau}
\displaystyle 
\vtau = \vx_s = \frac{\vx_\rho}{|\vx_\rho|}, \qquad \vnu = \vtau^\perp
\end{equation}
are a  unit tangent and unit normal to $\Gamma(t)$ respectively, 
where $(\cdot)^\perp$ denotes counter-clockwise rotation by $\frac{\pi}{2}$. 
The normal and tangential velocities of $\Gamma(t)$ are then
\begin{equation} \label{vpsi}
v= \vx_t \cdot \vnu,\qquad\psi= \vx_t \cdot \vtau.
\end{equation}
Note that (\ref{Gamma1}) only prescribes
$v$ and with it the shape of $\Gamma(t)$, so that there is a certain freedom 
in choosing the tangential velocity. Since $\vx_{ss}=\kappa\, \vnu$ one may consider
\begin{equation}  \label{xnormal}
\displaystyle
\vx_t = \vx_{ss} + f(\widetilde{w})\,\vnu = \frac{1}{| \vx_{\rho} | } 
%\frac{d}{d \rho} 
\left( \frac{ \vx_{\rho}}{| \vx_{\rho} | } \right)_\rho + f(\widetilde{w})\, \vnu,
\end{equation}
where $\widetilde{w}(\rho,t):=w(\vx(\rho,t),t), (\rho,t) \in \Ip \times [0,T]$. 
Clearly, (\ref{Gamma1}) holds with
the additional property that the velocity vector $\vx_t$ points in the normal direction. 
Coupling (\ref{xnormal})
to the PDE satisfied by $\tilde{w}$, Pozzi and Stinner have derived and analysed in  
\cite{PozziS15} a finite element scheme for (\ref{Gamma1},b) (with $g \equiv 0$)
based on continuous piecewise linears.  They are able to
prove the following error bounds in the spatially discrete case:
\begin{displaymath}
\sup_{[0,T]} \int_{\Ip} \left( | \vx_{\rho} - \vx^h_{\rho} |^2  
+ | \widetilde{w} - \widetilde{w}^h |^2 \right) \drho
+ \int_0^T \int_{\Ip} \left(
| \vx_t - \vx^h_{t} |^2 + | \widetilde{w}_{\rho} - \widetilde{w}^h_{\rho}|^2 \right)
\drho \dt \leq C \,h^2.
\end{displaymath}
A major difficulty in the analysis arises from the fact that (\ref{xnormal}) is only 
weakly parabolic. Following \cite{Dziuk94},
this problem is solved in \cite{PozziS15}  by deriving additional equations for the continuous 
and discrete length elements and by splitting the error $\vx_{\rho}-\vx^h_{\rho}$ into 
a tangent part and a length element part.
Apart from these analytical difficulties, the motion in purely the normal direction also 
may lead to the accumulation of
mesh points in numerical simulations. 
A natural way to handle the above mentioned difficulties is to introduce a tangential
part in the velocity, which can be seen as a reparametrisation, an approach that has recently been 
explored in a systematic way by Elliott and Fritz
in \cite{ElliottF15}. The underlying idea uses the DeTurck trick 
in coupling the motion of the curve to the harmonic
map heat flow, which results in the following equation replacing (\ref{xnormal}) (cf. (3.1) in \cite{ElliottF15}):
\begin{equation}  \label{xtang}
\displaystyle
\alpha \,\vx_t + (1-\alpha)\,(\vx_t \cdot \vnu)\, \vnu = \frac{\vx_{\rho \rho}}{| \vx_{\rho} |^2} 
+ f(\widetilde{w}) \,\vnu.
\end{equation}
In the above, 
$\alpha \in (0,1]$ is a parameter so 
that $\frac{1}{\alpha}$ corresponds to the diffusion coefficient in the harmonic map heat flow. 
Note that we
obtain (\ref{Gamma1}) by taking the scalar product of (\ref{xtang}) with $\vnu$. 
It turns out that $\vx(\cdot,t)$ gets
closer to a parametrisation proportional to arclength as $\alpha$ gets small. 
At the numerical level this means that mesh points along
the curve become more and more equidistributed. Setting formally $\alpha = 0$ one recovers an
approach introduced by Barrett, Garcke and N\"urnberg, see \cite{triplej}, \cite{fdfi}. 
A nice feature of (\ref{xtang}) is that
the problem now is strictly parabolic allowing for a more straightforward error analysis;
for the spatially discrete case with $f=0$, 
 see \cite{DeckelnickD95}, \cite{DeckelnickE98} for $\alpha=1$,  
and \cite{ElliottF15} for $\alpha \in (0,1]$. 

It remains to derive the equation for $\widetilde{w}$, which in contrast to \cite{PozziS15} 
will also involve the tangential velocity $(\vx_t \cdot\vtau)$. It is easily seen that
\begin{displaymath}
\widetilde{w}_t(\rho,t) = \frac{d}{dt}[ w(\vx(\rho,t),t)] = \matpartv w(\vx(\rho,t),t)  
+ ( \vx_t(\rho,t) \cdot \vtau(\vx(\rho,t)) \, w_s(\vx(\rho,t),t).
\end{displaymath}
Inserting this relation into (\ref{Gamma2}) and recalling (\ref{vpsi}), we obtain
\begin{equation}  \label{tildew}
\widetilde{w}_t - (\vx_t \cdot \vtau)\, 
\frac{1}{| \vx_{\rho} |} \,\widetilde{w}_{\rho} - \frac{d}{| \vx_{\rho} |} 
%\frac{\partial}{\partial \rho} 
\left( \frac{\widetilde{w}_{\rho}}{| \vx_{\rho} |} \right)_\rho 
- \kappa \, v \, \widetilde{w} = g(v,\widetilde{w}).
\end{equation}
The initial conditions for (\ref{xtang}), (\ref{tildew}) now are: $\vx(\rho,0)=\vx^0(\rho)$, 
$\widetilde{w}(\rho,0)=\widetilde{w}^0(\rho):=w^0(\vec{x}^0(\rho))$, $\rho \in \Ip$, 
where $\vec{x}^0$ is
a parametrization of $\Gamma^0$. 
We shall derive in Section \ref{sec:2} a weak formulation for (\ref{xtang}), 
(\ref{tildew}) which forms the basis for a discretisation
by continuous piecewise linear finite elements in space and 
a backward Euler scheme in time. As the main result of our paper we shall
present an error analysis for the resulting fully discrete scheme. 
The precise result will be formulated at the end of Section \ref{sec:2},
while the proof is carried out in detail in Section \ref{sec:3}. 
In Section \ref{sec:5} we present a series of test calculations that
confirm our theoretical results and demonstrate the above mentioned improvement of the mesh quality 
for small values of $\alpha$. 

Finally, we end this section with a few comments about notation.
We adopt the standard notation 
for Sobolev spaces, denoting the norm of
$W^{\ell,p}(G)$ ($\ell \in {\mathbb N}$, $p \in [1, \infty]$ 
and $G$ a bounded 
interval in ${\mathbb R}$) by
$\|\cdot \|_{\ell,p,G}$ and the semi-norm by $|\cdot |_{\ell,p,G}$. For
$p=2$, $W^{\ell,2}(G)$ will be denoted by
$H^{\ell}(G)$ with the associated norm and semi-norm written,
as respectively,  
$\|\cdot\|_{\ell,G}$ and $|\cdot|_{\ell,G}$.
For ease of notation, in the common case when $G\equiv \Ip$
the subscript ``$\Ip$'' will be dropped on the above norms and semi-norms.
%We introduce also $(H^{1}(\Omega))'$,
%as the dual space of $H^{1}(\Omega)$,
%and denote its norm by $\|\cdot\|_{(H^1)'}$.
The above are naturally extended to vector functions, and we will write $[W^{\ell,p}(G)]^2$
for a vector function with two components.
%and
%$\langle \cdot, \cdot \rangle$ denotes the duality pairing between $(H^1(\Omega))'$ and
%$H^1(\Omega)$. 
In addition, we adopt the standard notation $W^{\ell,p}(a,b;X)$
($\ell \in {\mathbb N}$, $p \in [1, \infty]$, $(a,b)$ an 
interval in ${\mathbb R}$, $X$ a Banach space) 
for time dependent spaces
with norm $\|\cdot\|_{W^{\ell,p}(a,b;X)}$.
Once again, we write $H^{\ell}(a,b;X)$ if $p=2$. 
Furthermore, $C$ denotes a generic constant independent of 
the mesh parameter $h$ and the time step $\Delta t$, see below.
%Finally,
%$C_{\gamma_1,\cdots,\gamma_\ell}$ denotes a constant dependent on 
%$\{\gamma_i^{-1}\}_{i=1}^{\ell}$.

\setcounter{equation}{0}
\section{Weak formulation and finite element approximation} \label{sec:2}
We shall assume that the data and the solution of (\ref{xtang}), (\ref{tildew}) (writing again $w$ instead
of $\widetilde{w}$ for ease of notation) are such that
\begin{subequations}
\begin{align}
&f \in C^1(\bR,\bR),\quad g\in C^1(\bR \times \bR,\bR),
\label{fgreg}\\
& \forall \ L >0, \ \exists \ C_L \geq 0 : 
\ | g(v_1,w) - g(v_2,w) | \leq C_L \,| v_1 - v_2 |
\quad \forall  \ v_1,\,v_2, \ |w| \leq L, 
 \label{glip} \\
&\vx \in W^{1,\infty}(0,T;[H^2(\Ip)]^2) %\cap W^{1,\infty}(0,T;[H^1(\Ip)]^2)
\cap H^2(0,T;[H^1(\Ip)]^2)\cap W^{2,\infty}(0,T;[L^2(\Ip)]^2),
\label{vxreg}\\
&w \in C([0,T];H^2(\Ip)) \cap W^{1,\infty}(0,T;H^{1}(\Ip))
\cap H^2(0,T;L^2(\Ip)), %\cap W^{1,\infty}(0,T;L^\infty(\Ip)),
\label{wreg} \\
&0 <m \leq |\vx_\rho| \leq M  \qquad \mbox{on }\Ip \times [0,T]
\label{nodegen}
\end{align}
for some $m,\,M \in \bR_{>0}$.
\end{subequations}
Since our error analysis will also include the discretisation in time our regularity assumptions are
slightly stronger than those made in \cite{PozziS15}, see Assumption 2.2 there.  
%However, they do not assume regularity on the second derivatives in time, i.e.\  
%$\vx \in H^{2}(0,T;[H^2(\Ip)]^2)$
%W^{1,\infty}(0,T;$ $[H^2(\Ip)]^2)$
%and $w \in H^2(0,T;L^2(\Ip))$;
%but they only prove an error bound for a semidiscrete scheme with no
%discretization in time. Hence, it is natural that we require this additional regularity.
%L^\infty(0,T;H^2(\Ip))\cap W^{1,\infty}(0,T;H^1(\Ip))$.
Let us fix $\alpha \in (0,1]$. 
For a test function 
$\vxi \in [H^1(\Ip)]^2$, we obtain on 
multiplying (\ref{xtang}) by $|\vx_\rho|^2\,\vxi$ that
\begin{equation}  \label{G1wb}
\displaystyle
\int_{\Ip} |\vx_\rho|^2 \left[\alpha \,\vx_t + (1-\alpha) 
\left(\vx_t \cdot \vnu \right) \vnu \right]\cdot \vxi \, \drho
+ \int_{\Ip} \vx_\rho \cdot \vxi_\rho \drho
= \int_{\Ip} |\vx_\rho|^2\,f(w)  \,\vnu \cdot \vxi\,\drho.
\end{equation} 
In order to derive a weak formulation for (\ref{tildew}) 
we employ an idea from \cite{DziukE07}, 
\cite{ElliottS12},  and calculate for $\eta \in H^1(\Ip)$ with the help of integration by parts and (\ref{tau})
\begin{eqnarray}
\lefteqn{
 \frac{d}{dt} \int_{\Ip} |\vx_\rho|\, w\,\eta \drho = \int_{\Ip} | \vx_\rho| \, w_t \, \eta \drho + 
\int_{\Ip} \frac{\vx_\rho \cdot \vx_{\rho,t}}{| \vx_\rho |} \, w \, \eta \drho }  \nonumber \\
&  = &  \int_{\Ip} | \vx_\rho| \, w_t \, \eta \, \drho - \int_{\Ip} 
( \vx_t \cdot \vtau) \,w_\rho \,\eta \drho - 
\int_{\Ip} (\vx_t \cdot \vtau) \, w \, \eta_\rho \drho - 
\int_{\Ip} \kappa \,( \vx_t \cdot \vnu)\, w \, \eta \, | \vx_\rho | \drho \nonumber \\
& & \qquad  =  -d \int_{\Ip} \frac{w_\rho\,\eta_\rho}{|\vx_\rho|} \drho 
- \int_{\Ip} %(\vx_t \cdot \vtau) 
\psi \,w\,\eta_\rho \drho
+  \int_{\Ip} |\vx_\rho|\,g\left(v,w\right) \eta \drho, \label{G2wb}
\end{eqnarray}
where we have also used (\ref{vpsi}), (\ref{tildew}) and the fact that
$\vtau_{\rho} %\cdot \vx_t 
= \kappa \,\vnu \,|\vx_\rho|$.

We now use (\ref{G1wb}), (\ref{G2wb}) in order to discretise our system and begin by
introducing the decomposition
$\Ip=\cup^J_{j=1}\overline{\sigma_j}$, where
$\sigma_j=(\rho_{j-1},\rho_j)$. We set 
$h := \max_{j=1,\ldots,J} h_j$,
where $h_j:=\rho_j-\rho_{j-1}$ and assume that
\begin{align}
h \leq C\,h_j, \qquad j=1,\ldots,J.
\label{hqu}
\end{align}
Let 
\begin{equation} \label{Vh1}
\Vh_1 := \{\chi \in C(\Ip) : \chi\!\mid_{\sigma_j} 
\mbox{ is affine},\ j=1,\ldots, J\} \subset H^1(\Ip)
\end{equation}
%\Vh_0 &:= \{\chi \in L^\infty(\Ip) : \chi\!\mid_{\sigma_j} 
%\mbox{ is constant},\ j=1\to J\},
%\label{Vh0}\\
%\Vh_{-1} &:= \{\chi \in L^\infty(\Ip) : \chi\!\mid_{\sigma_j} 
%\mbox{ is affine},\  j=1\to J\}.
%\label{Vh-1}
%\end{align}
%\end{subequations}
and $I^h:C(\Ip) \to \Vh_1$  be the
standard Lagrange interpolation operator such that $(I^h \eta)(\rho_j)=\eta(\rho_j)$,
$j=1,\ldots,J$. We require also the local interpolation operator $I^h_j \equiv I^h_{| \sigma_j}$,
$j=1, \ldots,J$ and recall for $p \in (1,\infty]$, $k \in \{0,1\}$,  
$\ell \in \{1,2\}$
and $j=1,\ldots,J$ that
\begin{subequations}
\begin{align}%{2}
h_j^{\frac{1}{p}}\,|\eta^h|_{0,\infty,\sigma_j} +
h_j\,|\eta^h|_{1,p,\sigma_j} &\leq C\,|\eta^h|_{0,p,\sigma_j} 
\qquad 
&&\forall \ \eta^h \in V^h_1,  \label{invh}\\
|(I-I^h_j)\eta|_{k,p,\sigma_j} &\leq C\,h^{\ell-k}_j\, |\eta|_{\ell,p,\sigma_j}
\qquad &&\forall \ \eta \in W^{\ell,p}(\sigma_j),
\label{Ih}\\
|(I-I^h_j)\eta|_{\ell-1,\infty,\sigma_j} &\leq C\,h^{\frac{1}{2}}_j\, |\eta|_{\ell,\sigma_j}
\qquad &&\forall \ \eta \in H^{\ell}(\sigma_j).
\label{Ihinf}
\end{align}
\end{subequations}  
As well as the standard $L^2(\Ip)$ inner product  $(\cdot,\cdot)$, 
we introduce the discrete inner product $(\cdot,\cdot)^h$ defined by
\begin{align}
\left( \eta_1, \eta_2 \right)^h :=
\sum_{j=1}^J \int_{\sigma_j} I^h_j(\eta_1\,\eta_2),
\label{innerh}
\end{align} 
where $\eta_i$ are piecewise continuous functions on the partition 
$\cup^J_{j=1}\overline{\sigma_j}$ of $\Ip$.
We note for $j=1,\ldots,J$ and for all $\eta^h, \chi^h \in V^h_1$ that
\begin{subequations}
\begin{align}
\int_{\sigma_j} |\eta^h|^2 \,\drho &\leq \int_{\sigma_j} I^h_j\left[|\eta^h|^2\right] \,\drho 
\leq  3\int_{\sigma_j} |\eta^h|^2 \,\drho, &&
\label{lumpev} \\
\left|\int_{\sigma_j} (I-I^h_j)(\eta^h\,\chi^h)\,\drho\right|
&\leq C\,h_j^2 \,|\eta^h|_{1,\sigma_j}\,|\chi^h|_{1,\sigma_j} \leq C\,h_j\,
|\eta^h|_{1,\sigma_j}\,|\chi^h|_{0,\sigma_j}.
\label{lumperr}
\end{align}
\end{subequations}
The result (\ref{lumperr}) follows immediately from (\ref{invh},b).
The inner products $(\cdot,\cdot)$ and $(\cdot,\cdot)^h$ are naturally extended to 
vector functions.   
In addition to the above spatial discretisation, let $0\equiv t_0<t_1<
\cdots<t_{N-1}<t_N \equiv T$ be a partitioning of $[0,T]$
with time steps $\Delta t_n:= t_n-t_{n-1}$,
$n=1,\ldots,N$, and $\Delta t := \max_{n=1,\ldots,N} \Delta t_n$. 
Before
we define our scheme we assign to an
element $\vX^n \in [V^h_1]^2$ (the upper index referring to the time level $n$) 
a piecewise constant discrete unit tangent and normal by
\begin{equation} \label{disctn}
\displaystyle 
\vTau^n = \frac{\vX_\rho^n}{|\vX_\rho^n|}, \qquad 
\vNu^n=(\vTau^n)^\perp \quad \mbox{ on }  \sigma_j,\; j=1,\ldots,J.
\end{equation}

Our discretisation of (\ref{G1wb}) now reads: given $\vX^{n-1} \in [V^h_1]^2$ and $W^{n-1} \in V^h_1$,
find $\vX^n \in [V^h_1]^2$ such that 
\begin{align}
&\left( |\vX^{n-1}_\rho|^2 \left[\alpha \,D_t \vX^n + (1-\alpha) 
\left( D_t \vX^n \cdot \vNu^{n-1} \right) 
\vNu^{n-1} \right] , \vxi^h \right)^h
+  \left( \vX_\rho^n , \vxi_\rho^h \right)  \nonumber \\
& \hspace{1in} =  \left( |\vX^{n-1}_\rho|^2\,f(W^{n-1})  
\,\vNu^{n-1} , \vxi^h\right)^h
\qquad \forall \ \vxi^h \in [V^h_1]^2.
\label{G1wbn}
\end{align}
Here, and in what follows,  we abbreviate $D_t a^n:= \frac{a^n - a^{n-1}}{\Delta t_n}$.
We next use the solution $\vX^n$ of (\ref{G1wbn}) in order to discretise (\ref{G2wb}). To do so,
we define approximations $V^n$, $\Psi^n$ of the normal and tangential velocities by
\begin{equation}  \label{vpsin}
\displaystyle 
V^n = D_t \vX^n \cdot \vNu^n
\quad \mbox{and}
\quad \Psi^n= D_t \vX^n \cdot \vTau^n \quad 
\mbox{ on } \sigma_j, \;  j=1,\ldots,J.
\end{equation}
Then find $W^{n} \in V^h_1$ such that 
\begin{align}
&D_t \left[\left( |\vX_\rho^n|\, W^n,\eta^h \right)^h \right] +  d\left(
\frac{W_\rho^n}{|\vX_\rho^n|},\eta^h_\rho\right)
+  \left(\Psi^n\,W^n,\eta_\rho^h \right)^h  \nonumber \\
&\hspace{1in} =  \left( |\vX_\rho^n|\,g(V^n,W^{n-1}),
\eta^h \right)^h \qquad \forall \ \eta^h \in V^h_1. \label{G2wbn} 
\end{align}

\noindent
Let us formulate
the main result of this paper, which will be proved in Section 3.
\begin{theorem} \label{main}
Let $\vX^0 = I^h \vx^0\in [V^h_1]^2$ and $W^0=I^h w^0 \in V^h_1$.
There exists $h^\star>0$ such that for all $h \in (0,h^\star]$ and $\Delta t \leq C\,h$ 
the discrete problem (\ref{G1wbn}), (\ref{G2wbn})
has a unique solution $(\vX^n,W^n) \in [V^h_1]^2 \times V^h_1, n=1,\ldots,N$, 
and the following error bounds hold:
\begin{align}
&\sup_{n=0,\ldots,N} \left[|\vx^n-\vX^n|_1^2 + |w^n-W^n|_0^2\right] 
\nonumber
\\
& \hspace{1in} + \sum_{n=1}^N \Delta t_n \left[ \left| \vx^n_t- D_t \vX^n \right|_0^2 
+ |w^n-W^n|_1^2 \right]
\leq C\,h^2,
\label{mainerr} 
\end{align} 
where $\vx^n:= \vx(\cdot,t_n),\, w^n:=w(\cdot,t_n),\, \vx_t^n:=\vx_t(\cdot,t_n), 
\,n=0,\ldots,N$.
\end{theorem}

\setcounter{equation}{0}
\section{Error analysis} \label{sec:3}
To begin, it follows from (\ref{vxreg},d) and (\ref{vpsi}) that for $n=0,\ldots,N$
\begin{align}
\|\vx^n\|_2  + \|\vx_t^n\|_1  + \|v^n\|_1 + \|\psi^n\|_1 +\|w^n\|_2 \leq C,
\label{vxwnreg}
\end{align} 
where $v^n=v(\cdot,t_n), \,\psi^n=\psi(\cdot,t_n)$. We abbreviate for $n=0,\ldots,N$ 
\begin{align}
\vE^n := I^h \vx^n - \vX^n \in [V^h_1]^2 \qquad
\mbox{and} \qquad Z^n := I^h w^n - W^n \in V^h_1.
\label{EnZn}
\end{align}

Our proof of Theorem \ref{main} will be based on induction.
We  assume for 
an $n \in \{1,\ldots,N\}$ that 
\begin{align}
%|\vE^{n-1}|_{1}^2 & \leq   h^2   \, e^{\gamma \int_0^{t_{n-1}} \zeta(s)ds} , 
%\label{Errind1} \\
|\vE^{n-1}|_{1}^2 + \beta^2 \,(| \vX_\rho^{n-1}\, | Z^{n-1}, Z^{n-1})^h & \leq   h^2   
 \, e^{\gamma \int_0^{t_{n-1}} \zeta(t)\dt}, \qquad 
 \mbox{for }\ h \in (0,h^{\star}],\label{Errind2}
\end{align}
where the function $\zeta$ is defined by
\begin{align}
\zeta(t) :=  1 +  \|\vx_{tt}(t)\|_1^2 +  | w_{tt}(t) |_0^2 
\label{zeta}
\end{align}
and $h^{\star}>0$ is chosen so small that
\begin{align}
&  (h^\star)^{\frac{1}{2}} \,e^{\gamma K} \leq 
\beta \quad \mbox{with}\quad K:= \int_0^T \zeta(t) \dt. \label{hstar}
\end{align}
Here,  $\beta \in (0,1]$ and $\gamma > 0$ are independent of $h$ and $\Delta t$, 
and will be chosen a posteriori.
Clearly, (\ref{Errind2}) holds for $n=1$ in view of our choice of initial data
for $\vX^0$ and $W^0$. \\
It follows from  (\ref{hqu}), (\ref{Ihinf}),  
(\ref{Errind2}), (\ref{vxwnreg}) and (\ref{hstar})  that 
\begin{align}
|\vx^{n-1}-\vX^{n-1}|_{1,\infty} &\leq |\vE^{n-1}|_{1,\infty} + |(I-I^h)\vx^{n-1}|_{1,\infty}
\leq C\,h^{-\frac{1}{2}}\,|\vE^{n-1}|_1 + C\,h^{\frac{1}{2}}\,|\vx^{n-1}|_{2}
\nonumber \\
&\leq C h^{\frac{1}{2}} 
\left( e^{\frac{\gamma K}{2}} +1 \right)  
\leq C\, (h^{\star})^{\frac{1}{2}} e^{\gamma K} \leq C \,\beta \leq
\min\left\{\frac{m}{2},M\right\},
\label{En3}
\end{align}
provided that $\beta$ is chosen small enough. Combining this inequality with 
(\ref{nodegen}) we infer that 
\begin{align}
&0<\frac{m}{2} \leq |\vX^{n-1}_\rho| \leq 2\,M \qquad \mbox{ on } \Ip  \mbox{ for }  h \in (0,h^\star].
\label{disnodegen}
\end{align}
If we use (\ref{disnodegen}) and argue similarly as in (\ref{En3}) we further obtain
for any 
$h \in (0,h^\star]$
%$0 < h \leq  h^\star$ 
\begin{subequations}
\begin{align}
 |\vtau^{n-1}-\vTau^{n-1}|_0 + |\vnu^{n-1} -\vNu^{n-1}|_0  
&\leq   C \, | \vx^{n-1} - \vX^{n-1} |_1  
\leq C \left[ |\vE^{n-1}|_1 + h \right], 
\label{tauerrnm}\\
|\vtau^{n-1}-\vTau^{n-1}|_{0,\infty} + |\vnu^{n-1} -\vNu^{n-1}|_{0,\infty}  
&\leq  C \, |  \vx^{n-1} - \vX^{n-1} |_{1,\infty} \nonumber \\
& \leq    C \,h^{-\frac{1}{2}}\left[ |\vE^{n-1}|_1 + h \right]. 
\label{tauerr1nm}
\end{align}
\end{subequations}
In the same way as (\ref{En3}), 
we obtain from  (\ref{lumpev}), (\ref{disnodegen}), (\ref{Errind2})
%(\ref{vxwnreg}) 
and (\ref{hstar}) that
\begin{equation} \label{Z0}
| Z^{n-1} |_0^2   \leq   C \left( | \vX_\rho^{n-1} |\, Z^{n-1}, Z^{n-1} \right)^h 
\leq \frac{C}{\beta^2}\, h^2 \,e^{\gamma K} \leq \frac{C}{\beta^2} \,h \,
h^{\star} e^{2 \gamma K} \leq C\,h,
\end{equation}
which combined with  (\ref{invh},c), (\ref{hqu}) and (\ref{vxwnreg}) yields for $h \in (0,h^\star]$
\begin{align}
& |W^{n-1}|_{0,\infty} \leq | Z^{n-1} |_{0,\infty} + | I^h w^{n-1} |_{0,\infty} 
\leq C\, h^{-\frac{1}{2}} | Z^{n-1} |_0 + C \leq C.
\label{disupw}
\end{align}

\subsection{The curve equation}\label{curerr}
We assume throughout that $h \in (0,h^\star]$. Since (\ref{G1wbn}) forms a linear problem it is 
easily seen that (\ref{disnodegen}) implies the existence of 
a unique solution $\vX^n \in [V^h_1]^2$ to (\ref{G1wbn}).
We deduce from (\ref{G1wbn}) and (\ref{G1wb}) with $\vxi=\vxi^h=
\Delta t_n \,D_t \vE^n \in [V^h_1]^2$ that 
\begin{align}
LHS & :=   \Delta t_n \left( |\vX^{n-1}_\rho|^2 \left[\alpha \,D_t \vE^n + (1-\alpha) 
\left( D_t \vE^n \cdot \vNu^{n-1} \right) 
\vNu^{n-1} \right] , D_t\vE^n \right)^h \nonumber \\
& \qquad  + \Delta t_n \left( \vE_\rho^n , (D_t\vE^n)_\rho \right)  \nonumber \\
&   =   \Delta t_n
\left( |\vX^{n-1}_\rho|^2 \left[\alpha \,D_t \vx^n + (1-\alpha) 
\left( D_t \vx^n \cdot \vNu^{n-1} \right) 
\vNu^{n-1} \right] , D_t\vE^n \right)^h \nonumber  \\
&  \qquad
+ \Delta t_n \left( \vx_\rho^n , (D_t\vE^n)_\rho \right)
- \Delta t_n \left( |\vX^{n-1}_\rho|^2\,f(W^{n-1})  
\,\vNu^{n-1} , D_t\vE^n\right)^h \nonumber \\
&   = \Delta t_n \biggl[ 
\left( |\vX^{n-1}_\rho|^2 \left[\alpha \,D_t \vx^n + (1-\alpha) 
\left( D_t \vx^n \cdot \vNu^{n-1} \right) 
\vNu^{n-1} \right] , D_t\vE^n \right)^h \nonumber  \\
& \qquad
-%\Delta t_n 
\left( |\vx^n_\rho|^2 \left[\alpha \,\vx_t(t_n) + (1-\alpha) 
\left(\vx_t(t_n) \cdot \vnu^{n} \right) 
\vnu^{n} \right] , D_t\vE^n \right) \biggr] \nonumber  \\
& \qquad  
+ \Delta t_n \left[
\left( |\vx^n_\rho|^2\,f(w^n)  
\,\vnu^n,D_t\vE^n\right) 
-\left( |\vX^{n-1}_\rho|^2\,f(W^{n-1})  
\,\vNu^{n-1} , D_t\vE^n\right)^h \right] \nonumber \\
&  =:   A_1 + A_2. \label{En} 
\end{align}
Using (\ref{disnodegen}) and (\ref{lumpev}) 
we find, with the help of an elementary calculation,
that
\begin{equation} \label{lhs}
LHS \geq \Delta t_n \,\alpha \, \frac{m^2}{4} \left| D_t\vE^n\right|_0^2 
+ \frac{1}{2} \left[ |\vE^n|_1^2+ |\vE^n-\vE^{n-1}|_1^2- | \vE^{n-1} |_1^2 \right].
\end{equation}
Let us analyse the $A_1$ term defined in (\ref{En}) and note that
\begin{align}
A_1 &=  \Delta t_n
\biggl[
\left( |\vX^{n-1}_\rho|^2 \left[\alpha \,D_t(I^h\vx^n) + (1-\alpha) 
\left( D_t (I^h \vx^n) \cdot \vNu^{n-1} \right) 
\vNu^{n-1} \right] , D_t\vE^n \right)^h
\nonumber \\
& \qquad  -
\left( |\vX^{n-1}_\rho|^2 \left[\alpha \, D_t (I^h \vx^n) + (1-\alpha) 
\left( D_t (I^h \vx^n) \cdot \vNu^{n-1} \right) 
\vNu^{n-1} \right] , D_t\vE^n \right) \biggr]
\nonumber \\
& \qquad  + \Delta t_n  \biggl[ \alpha
\left( |\vX^{n-1}_\rho|^2 \left[ D_t (I^h \vx^n) - \vx_t^n\right] , D_t\vE^n \right)
\nonumber \\
& \qquad  
+ (1-\alpha) \left( |\vX^{n-1}_\rho|^2
\left( \left[ D_t (I^h \vx^n) - \vx_t^n\right] \cdot \vNu^{n-1} \right) 
\vNu^{n-1} , D_t\vE^n \right) \biggr]
\nonumber \\
& \qquad  + \Delta t_n \biggl[ (1-\alpha)\, \left( |\vX^{n-1}_\rho|^2\,
\left[ \left(\vx_t^n \cdot \vNu^{n-1} \right) 
\vNu^{n-1} - \left(\vx_t^n \cdot \vnu^{n}\right) 
\vnu^{n}\right], D_t\vE^n \right) \nonumber\\
& \qquad
+  
\left(\left[|\vX^{n-1}_\rho|^2 - |\vx^n_\rho|^2\right] \left[\alpha \,\vx_t^n + (1-\alpha) 
\left(\vx_t^n \cdot \vnu^{n} \right) 
\vnu^{n} \right] , D_t\vE^n \right) \biggr]
  =: \sum_{i=1}^3 A_{1,i}. 
\label{A1}
\end{align}
We now bound the terms $A_{1,i}$ defined in (\ref{A1}) on recalling
(\ref{vxreg},e), 
(\ref{disnodegen}), (\ref{tau}) and (\ref{disctn}).
It follows from  (\ref{lumperr}) and (\ref{Ih}) that
%\footnote{Be consistent in showing the dependence on the generic constant $C$
%on norms of $\vx$ and $w$. We are very precise early on, then not.}
\begin{align}
|A_{1,1}| &\leq C\,h\,\Delta t_n \,|D_t (I^h \vx^n)|_{1}\,|D_t\vE^n|_{0}
\leq C\,h\,\Delta t_n   
\,|D_t\vE^n|_{0},
\label{A11}
\end{align}
since $| D_t(I^h \vx^n) |_1 \leq | D_t \vx^n |_1 \leq \|\vx_t\|_{L^\infty(0,T;[H^1(\Ip)]^2)} \leq C$.
Similarly, 
\begin{align}
|A_{1,2}| & \leq C\, \Delta t_n  \left[
\,|(I-I^h) D_t \vx^n|_{0} 
+ | D_t \vx^n - \vx_t^n|_{0} \right]\,|D_t\vE^n|_{0}
\nonumber \\
&\leq C \,\Delta t_n  \left[ h %\, | \vx^n_t |_1 
+  \left(\Delta t_n \int_{t_{n-1}}^{t_n}|\vx_{tt}|_0^2 \dt\right)^{\frac{1}{2}}
\right]   
|D_t\vE^n|_{0}.
\label{A12}
\end{align}
Next, we have from Sobolev embedding, (\ref{Ih}) and (\ref{vxwnreg}) that
\begin{align} 
|A_{1,3}| & \leq  C \,\Delta t_n\, 
\left[ |\vnu^{n} - \vNu^{n-1}|_{0}
+ \left| |\vx^{n}_\rho|^2 - |\vX^{n-1}_\rho|^2 \right|_{0}
\right]\,|D_t\vE^n|_{0}   
\nonumber \\
&\leq   C \,\Delta t_n \left[ 
%| \vx^{n}-\vx^{n-1} |_1 
\Delta t_n\,|D_t \vx^n|_1 +
| \vx^{n-1}- \vX^{n-1}|_{1} \right] 
%\|\vx\|_{H^2(0,T;H^1(\Ip))}\,
|D_t\vE^n|_{0} \nonumber \\
& \leq   C  \,\Delta t_n \left[ 
h + | \vE^{n-1} |_1 \right] | D_t \vE^n |_0, 
\label{A134}
\end{align}
since $\Delta t_n \leq C h$.
We now turn our attention to the $A_2$ term in (\ref{En}) and note that 
\begin{align}
A_2 &=   \Delta t_n 
\left( |\vx^n_\rho|^2\,f(w^n)  
\,\left(\vnu^n-\vNu^{n-1}\right) + \left(|\vx^n_\rho|^2-|\vX^{n-1}_\rho|^2\right)f(w^n)  
\,\vNu^{n-1} ,D_t\vE^n\right) 
\nonumber \\
& \qquad + \Delta t_n \biggl[ \left( |\vX^{n-1}_\rho|^2\,f(W^{n-1})  
\,\vNu^{n-1} , D_t\vE^n\right)\nonumber \\
& \qquad \qquad \qquad \qquad -\left( |\vX^{n-1}_\rho|^2\,f(W^{n-1})  
\,\vNu^{n-1} , D_t\vE^n\right)^h \biggr]
\nonumber \\
&  \qquad + \Delta t_n
\left( |\vX^{n-1}_\rho|^2\left( f(w^n) -f(W^{n-1})\right)\,\vNu^{n-1},D_t\vE^n\right)
=:  \sum_{i=1}^3 A_{2,i}.
\label{A2} 
\end{align}
We bound the terms $A_{2,i}$ defined in (\ref{A2}) on recalling
(\ref{fgreg}--e),  
(\ref{disnodegen}), (\ref{disupw}), (\ref{tau}) and (\ref{disctn}).
Since $|f(w^n)|_{0,\infty} \leq C$, it follows  similarly to (\ref{A134}) that 
\begin{align}
|A_{2,1}| &\leq C\,\Delta t_n 
%_{L^\infty(t_{n-1},t_n;L^\infty(\Ip))}
\,|\vx^n-\vX^{n-1}|_1\,|D_t\vE^n|_0
\leq C \,\Delta t_n \left[ h +  |\vE^{n-1}|_{1} \right] |D_t\vE^n|_{0}.
\label{A212} 
\end{align}
Next, %similarly to (\ref{A11}), 
we have from 
(\ref{innerh}), (\ref{lumperr}), (\ref{invh},b), %(\ref{disupw}), 
(\ref{hqu}) and (\ref{vxwnreg}) that
\begin{align}
|A_{2,2}| &\leq C\,\Delta t_n
\left[\, |(I-I^h)f(W^{n-1})|_0 +
h\,| I^h[f(W^{n-1})]|_1\right]|D_t\vE^n|_0 \nonumber \\
& \leq C\,\Delta t_n\,h\,|f'(W^{n-1})|_{0,\infty}\,|W^{n-1}|_1\,|D_t\vE^n|_0
\nonumber \\ &
\leq C\,\Delta t_n\,h\left[ \,| I^h w^{n-1}|_1 + |Z^{n-1}|_1   \right]|D_t\vE^n|_0
\leq C\,\Delta t_n\left[\, h + |Z^{n-1}|_0 \right]|D_t\vE^n|_0.
\label{A24}
\end{align}
Finally, it follows from %(\ref{fgreg}), (\ref{disupw}), 
(\ref{Ih}), (\ref{vxwnreg}) and as $\Delta t_n \leq Ch$  that 
\begin{align} 
|A_{2,3}| &\leq C\,\Delta t_n\,|f(w^n)-f(W^{n-1})|_0\,|D_t\vE^n|_0 \leq 
C\,\Delta t_n\,|w^n-W^{n-1}|_0\,|D_t\vE^n|_0 
\nonumber \\
&\leq 
C\,\Delta t_n \left[ \Delta t_n  +  h \, | w^{n-1} |_1 + |Z^{n-1}|_0\right]\,|D_t\vE^n|_0 \nonumber \\
& \leq   C \,\Delta t_n \left[ h +  |Z^{n-1}|_{0}\right]|D_t\vE^n|_0.
\label{A23} 
\end{align}
If we combine (\ref{A1})--(\ref{A23}) with Young's inequality and the definition of $\zeta$ we infer,
as $\Delta t_n \leq C\,h$, that
\begin{align} 
| A_1 + A_2 | &\leq \delta \,\Delta t_n \,| D_t \vE^n |_0^2 
+ C(\delta) \left[ \Delta t_n \left(| \vE^{n-1} |_1^2 + | Z^{n-1} |_0^2\right) + 
h^2 \int_{t_{n-1}}^{t_n} \zeta(t) \dt \right].
\label{rhs}
\end{align}
Inserting (\ref{lhs}) and  (\ref{rhs}) into  (\ref{En}) and choosing $\delta$ sufficiently small
we obtain with the help of (\ref{disnodegen}) and (\ref{lumpev})
\begin{align}
&\Delta t_n \,\alpha \, \frac{m^2}{4} \left| %\frac{\vE^n-\vE^{n-1}}{\Delta t_n} 
D_t\vE^n\right|_0^2 
+ |\vE^n|_1^2+ |\vE^n-\vE^{n-1}|_1^2  \nonumber\\
& \; \leq  | \vE^{n-1}|_1^2 + 
C \, \Delta t_n  \left[  | \vE^{n-1}|_1^2 +   (| \vX_\rho^{n-1}|\, Z^{n-1}, Z^{n-1})^h   \right]  
+ C \, h^2 \int_{t_{n-1}}^{t_n} \zeta(t) \dt.
\label{En1}
\end{align}
The induction hypothesis (\ref{Errind2}) together with the fact that $\zeta \geq 1$ then yields
\begin{align} 
&\Delta t_n \,\alpha \, \frac{m^2}{4} \left| %\frac{\vE^n-\vE^{n-1}}{\Delta t_n} 
D_t\vE^n\right|_0^2 
+ |\vE^n|_1^2+ |\vE^n-\vE^{n-1}|_1^2 \nonumber \\
& \hspace{1in} \leq   h^2 \,   e^{\gamma \int_0^{t_{n-1}} \zeta(t) \dt} 
\left[ 1  + C \left( 1 + \frac{1}{\beta^2} \right)\int_{t_{n-1}}^{t_n} 
\zeta(t) \dt \right] \label{En4} 
\end{align}
so that 
\begin{align}
& | \vE^n |_1^2 \leq h^2 \,  e^{\gamma \int_0^{t_n} \zeta(t) \dt} 
\quad \mbox{ provided that } \quad \gamma \geq C \left( 1 + \frac{1}{\beta^2} \right). \label{gamma1}
\end{align}
In particular, (\ref{En4}), this choice of $\gamma$ and (\ref{hstar}) imply that
\begin{equation}  \label{Enh}
\displaystyle \Delta t_n \,\alpha \, \frac{m^2}{4} \left| %\frac{\vE^n-\vE^{n-1}}{\Delta t_n} 
D_t\vE^n\right|_0^2 
+ |\vE^n|_1^2+ |\vE^n-\vE^{n-1}|_1^2 \leq h^2 \,e^{\gamma K} \leq h^{\frac{3}{2}} \,
(h^{\star})^{\frac{1}{2}} \,e^{\gamma K} \leq \beta \,h^{\frac{3}{2}}.
\end{equation}

\begin{comment}
& \quad \leq 
C_1 \exp(C_1\,t_k)\sum_{n=1}^k \Delta t_n 
\,|Z^{n-1}|_0^2 
\nonumber \\
& \quad \quad +
C_2 \exp(C_1\,t_k)\left[h^2\,t_k + h^4 \int_{0}^{t_k} |\vx_t|^2_2 \dt
+ (\Delta t)^2 \int_{0}^{t_k} \left[
\,|\vx_{t}|_1^2
+ |\vx_{tt}|_0^2 + |w_{t}|_0^2 \, \right]\dt \right]
\nonumber \\
& \quad \leq 
\exp(C_1\,T) \left[C_3 + C_1 \sum_{n=1}^k \Delta t_n\,C^Z_{n-1}\right] 
\left(h^2 + (\Delta t)^2\right) \leq C^E_k h^2
\label{En2}
\end{align}
for some $C^E_k \in \bR_{>0}$, independent of $h$ and $\Delta t$.  
\end{comment}

In the same way as in (\ref{En3})--(\ref{tauerrnm},b) 
we infer from (\ref{gamma1}) for $h \in (0,h^\star]$
\begin{subequations}
\begin{align}
0<\frac{m}{2}  \leq |\vX^{n}_\rho| & \leq 2\,M \quad \mbox{ on } \Ip,
\label{disnodegenk}\\
|\vtau^n-\vTau^n|_0 
+  |\vnu^n -\vNu^n|_0 
+ h^{\frac{1}{2}} \left(|\vtau^n-\vTau^n|_{0,\infty}  
+ |\vnu^n -\vNu^n|_{0,\infty}\right) 
&\leq C \left[ |\vE^n|_1 + h \right].
\label{tauerr}
\end{align}
\end{subequations}
In addition, we deduce from (\ref{vpsi}), (\ref{vpsin}), (\ref{vxwnreg}),
(\ref{tauerr}), (\ref{vxreg}) and as $\Delta t_n \leq C\,h$ 
that
%\footnote{Need extra regularity than 
%(\ref{vxreg}) for this. Require in addition that $\vx \in W^{2,\infty}(0,T;[L^2(\Ip)]^2)$.}
\begin{align}
&|v^n-V^n|_0 + |\psi^n-\Psi^n|_0  \nonumber \\
& \quad \leq   
|\vx_t^n\cdot(\vtau^n-\vTau^n)|_0 + |\vx_t^n\cdot(\vnu^n-\vNu^n)|_0
 + 2 \left|\vx_t^n - D_t (I^h \vx^n) \right|_0 
 +2 \left|%\frac{\vE^n-\vE^{n-1}}{\Delta t_n}
 D_t\vE^n\right|_0
 \nonumber \\
& \quad \leq  
C \left[ | \vE^n |_1 + h \right]
%+ C\,h \, |\vx^n_{t} |_1
 + C\,\Delta t_n\,\|\vx_{tt}\|_{L^\infty(t_{n-1},t_n;[L^2(\Ip)]^2)}+  
2\left|
D_t\vE^n\right|_0
\nonumber \\
& \quad \leq 
C \left[ | \vE^n |_1 + h +  
\left| D_t\vE^n\right|_0 \right].
\label{vpsierr}
\end{align}
Furthermore, we conclude from  (\ref{invh},b), (\ref{hqu}), (\ref{vpsierr}),  
(\ref{vxwnreg}), (\ref{Enh}) and as $\Delta t_n \leq C\,h$ that
\begin{align}
\Delta t_n \, |\Psi^n|^2_{0,\infty} &\leq 2\,\Delta t_n 
\left[|I^h\psi^n|^2_{0,\infty}+|I^h\psi^n-\Psi^n|^2_{0,\infty}
\right] \nonumber \\
& \leq  2\,\Delta t_n \left[|\psi^n|^2_{0,\infty}+ C\,h^{-1}\,|I^h\psi^n-\Psi^n|^2_{0} \right]
\nonumber \\
&\leq 2\,\Delta t_n \,|\psi^n|^2_{0,\infty} 
+  C\,\Delta t_n\,h^{-1} \left[|\psi^n-\Psi^n|^2_{0} + |(I-I^h)\psi^n|_0^2 \right]
\nonumber \\
&\leq C \,h\, \Vert \psi^n \Vert^2_1
+  C \, \Delta t_n \, h^{-1} \left[  | \vE^n |^2_1 + h^2 \right]
+ C \,h^{-1} \Delta t_n \left|D_t\vE^n \right|^2_0 \nonumber \\
& \leq C\,  h^{\frac{1}{2}}  < 4\,d 
\label{Psinbd}
\end{align}
for $h \in (0, h^\star]$, provided that $h^\star$ is  small enough. 
Finally, similarly to (\ref{A11}), we have on noting (\ref{invh},b), (\ref{hqu}) and (\ref{vxreg}) 
that
\begin{equation}  \label{DtX}
\displaystyle | \vX^n - \vX^{n-1} |_1 \leq C \,\Delta t_n 
\left[ | D_t \vE^n |_1 + | D_t (I_h \vx^n) |_1 \right] \leq 
C \,\Delta t_n \left[ h^{-1} \,| D_t \vE^n |_0 + 1\right].
\end{equation}
 
\begin{comment}  
\begin{rem}
We remark that the assumption (\ref{Deltath1}) is required to obtain the second bound in 
(\ref{Deltabdk}), which is required if $g$ depends on $v$; see (\ref{B4}) below.
If $g$ depends solely on $w$, then the assumption (\ref{Deltath1}) is not required.
\end{rem}
\end{comment}

\subsection{The scalar equation}\label{scaerr}
We assume throughout that $h \in (0,h^\star]$. 
Let us first establish the existence and uniqueness of 
$W^n \in V^h_1$. Since the system (\ref{G2wbn}) is
linear, existence follows from uniqueness which in turn is a consequence of the estimate
\begin{align*}
 \left|\left(\Psi^n\,\eta^h,\eta^h_\rho\right)^h\right| & \leq  
d\left(\frac{\eta^h_\rho}{|\vX_\rho^n|},\eta^h_\rho\right)
+\frac{1}{4\,d}\left( |\vX_\rho^n|\,(\Psi^n)^2\, \eta^h,\eta^h \right)^h  \\
& <  d\left(\frac{\eta^h_\rho}{|\vX_\rho^n|},
\eta^h_\rho\right) + \frac{1}{\Delta t_n} \left( | \vX^n_\rho | \, \eta^h,\eta^h \right)^h
\end{align*}  
in view of (\ref{Psinbd}). 
We deduce from (\ref{G2wbn}) 
with $\eta^h = Z^n$ 
that
%for all $Z^n \in V^h_1$ 
\begin{align}
&\left( D_t \left[ |\vX_\rho^n|\, Z^n \right],Z^n \right)^h + d\left(
\frac{Z_\rho^n}{|\vX_\rho^n|},Z^n_\rho\right)
 = \left( D_t \left[ |\vX_\rho^n|\, I^h w^n \right],Z^n \right)^h \nonumber \\
& \hspace{0.4in} + d \left( \frac{(I^h w^n)_\rho}{|\vX_\rho^n|},Z^n_\rho\right) 
+  \left(\Psi^n\, W^n,Z^n_\rho \right)^h
- \left( |\vX_\rho^n|\,g(V^n,W^{n-1}),Z^n\right)^h. \label{Zn0}
\end{align}
Observing that for all $a^{n-1},\,a^n,\,b^n,\,b^{n-1} \in \bR$ 
\begin{displaymath}
D_t [a^n\,b^n]\,b^n = 
\tfrac{1}{2} D_t [ a^n\,(b^n)^2]
+ \tfrac{1}{2 \Delta t_n }\,a^{n-1}\,(b^n-b^{n-1})^2 +
\tfrac{1}{2 \Delta t_n}\,(a^n-a^{n-1})\,(b^{n})^2
\end{displaymath}
and multiplying (\ref{Zn0}) by $\Delta t_n$ we obtain with the help of (\ref{G2wb}) and (\ref{disnodegenk})
\begin{eqnarray}
\lefteqn{ \tfrac{1}{2}\,  (|\vX^n_\rho|\,Z^n,Z^n)^h 
+  \frac{d}{2M} \, \Delta t_n  \, | Z^n |_1^2 
+ \tfrac{1}{2} \,(|\vX^{n-1}_\rho|\,(Z^n-Z^{n-1}),(Z^n-Z^{n-1})\,)^h } \nonumber \\
& \leq & \tfrac{1}{2} \, (|\vX^{n-1}_\rho|\,Z^{n-1},Z^{n-1})^h
- \tfrac{1}{2}\,\left(\left[|\vX^n_\rho|-|\vX^{n-1}_\rho|\right]Z^n,Z^n\right)^h
\nonumber \\
& & \quad   +   \left[\Delta t_n \, \left( D_t [|\vX_\rho^n|\, I^h w^n],Z^n \right)^h
- \Delta t_n \left( \left(|\vx_\rho|\,w\right)_t(t_n),Z^n\right)\right]
\nonumber \\
& & \quad  + \Delta t_n\,d \left(
\frac{(I^h w^n)_\rho}{|\vX_\rho^n|}-\frac{\,w_\rho^n}{|\vx_\rho^n|},Z^n_\rho\right)
%\nonumber \\
%& \quad \qquad 
+ \Delta t_n \left[\left(\Psi^n\,(I^h  w^n ),Z^n_\rho \right)^h
- \left(\psi^n\,w^n,Z^n_\rho \right)\right]
\nonumber \\
& &  \quad  - \Delta t_n \left(\Psi^n\,Z^n,Z^n_\rho \right)^h + \Delta t_n\left[ \left(|\vx_\rho^n|\,g(v^n,w^n),Z^n\right)
-\left( |\vX_\rho^n|\,g(V^n,W^{n-1}),Z^n\right)^h\right] 
\nonumber \\
& & \quad \qquad =:  \tfrac{1}{2}\,  (|\vX^{n-1}_\rho|\,Z^{n-1},Z^{n-1})^h +  \sum_{1=1}^6 B_i.
\label{Zn}
\end{eqnarray}

We now estimate the $B_i$ terms defined in (\ref{Zn}). In order to treat $B_1$ 
we first observe that
\begin{align*}
| \vX^n_\rho | - | \vX^{n-1}_\rho | = \vTau^{n-1} \cdot (\vX^n_\rho - \vX^{n-1}_\rho ) + 
\tfrac{1}{2} \,| \vX^n_\rho | \, | \vTau^n - \vTau^{n-1} |^2,
\end{align*}
so that (\ref{lumperr}), (\ref{invh}), (\ref{hqu}) and integration by parts  imply
\begin{eqnarray}
2 B_1 & \leq &  - \left( \vTau^{n-1} \cdot (\vX^n_\rho - \vX^{n-1}_\rho ) \,Z^n, Z^n \right)^h 
  \nonumber  \\
& \leq &   - \left( \vTau^{n-1} \cdot (\vX^n_\rho - 
\vX^{n-1}_\rho ) \,Z^n, Z^n \right) + C\, h\, | \vX^n - \vX^{n-1} |_1 \, | Z^n |_{0,\infty} \,  
| Z^n |_1 
\nonumber \\
& = & \left( (  \vtau^{n-1} - \vTau^{n-1})\cdot (\vX^n_\rho - \vX^{n-1}_\rho )\, Z^n, Z^n \right) 
 + \left( \vtau^{n-1}_\rho \cdot (\vX^n - \vX^{n-1})\, Z^n, Z^n \right) \nonumber \\
 & & +2  \left( \vtau^{n-1} \cdot (\vX^n - \vX^{n-1} )\, Z^n, Z^n_\rho \right) +
 C\, h\, | \vX^n - \vX^{n-1} |_1 \, | Z^n |_{0,\infty} \,  
| Z^n |_1.
\nonumber
\end{eqnarray}
Using (\ref{invh}), (\ref{hqu}), (\ref{tauerrnm}), (\ref{DtX}), 
the inequality $| \eta |_{0,\infty}^2 \leq C \,| \eta |_0\, 
\Vert \eta \Vert_1$, (\ref{A11})  and the fact that  
$| \vE^{n-1} |_1 \leq  h^{\frac{3}{4}}$ (cf.\ (\ref{Errind2}) and (\ref{hstar}))
we infer that
\begin{align}
B_1  & \leq C\left[|\vtau^{n-1}-\vTau^{n-1}|_0\,|\vX^n-\vX^{n-1}|_{1}\,|Z^n|_{0,\infty}+ 
|\vX^n-\vX^{n-1}|_0\, \Vert Z^n \Vert_1\right]|Z^n|_{0,\infty}\nonumber   \\
& \leq C \,\Delta t_n \left[ | \vE^{n-1} |_{1}+h \right] \left[ h^{-1} 
| D_t \vE^n |_0 +1 \right] | Z^n |_0 \Vert Z^n \Vert_1 \nonumber \\
& \qquad + C \,\Delta t_n \,|D_t(I^h \vx^n)|_0\,
| Z^n |_0^{\frac{1}{2}} \,\Vert Z^n \Vert_1^{\frac{3}{2}} 
+C \,\Delta t_n  \,| D_t \vE^n |_0 \, h^{-\frac{1}{2}}\,
| Z^n |_0 \Vert Z^n \Vert_1   \nonumber \\
& \leq C \,\Delta t_n \left[ h^{-\frac{1}{2}} \,| D_t \vE^n |_0 + h^{\frac{3}{4}} \right] 
| Z^n |_0 \, \Vert Z^n \Vert_1 + C \,\Delta t_n \,| Z^n |_0^{\frac{1}{2}} \,
\Vert Z^n \Vert_1^{\frac{3}{2}}.
\nonumber \\
\label{B1} 
\end{align}
Let us postpone the rather complicated analysis of $B_2$ and first deal with $B_3,\ldots,B_6$.
It follows from (\ref{vxreg}--e), (\ref{disnodegenk}), (\ref{Ih}) and Sobolev embedding  
that
\begin{align}
B_3 &=
\Delta t_n\,d\left(
\frac{\left[|\vx^n_\rho|- |\vX^n_\rho|\right]}{|\vx^n_\rho|\,|\vX^n_\rho|}w^n_\rho,  
Z^n_\rho
\right)
\leq C\,\Delta t_n\,
|\vx^n-\vX^n|_1 \,|w^n|_{1,\infty}\,|Z^n|_1
\nonumber \\
&\leq C\,\Delta t_n
\left[ 
|\vE^n|_1 + |(I-I^h)\vx^n|_{1}\right]|Z^n|_1 
\leq C\,\Delta t_n
\left[ %h\,|\vx^n|_2 +
|\vE^n|_1 + h\right]|Z^n|_1.
\label{B3}
\end{align}
Next we note from  (\ref{lumpev},b), (\ref{invh},b), (\ref{hqu}),
Sobolev embedding, (\ref{vxwnreg}) and (\ref{vpsierr}) that
\begin{align}
B_4 &=   \Delta t_n \left((\Psi^n-I^h \psi^n)\,I^h w^n,Z^n_\rho\right)^h 
 \nonumber \\
&   \qquad + \Delta t_n \, \left[\left((I^h \psi^n)\,(I^h w^n ),Z^n_\rho \right)^h
- \left((I^h \psi^n)\,(I^h w^n),Z^n_\rho \right) \right] 
\nonumber \\
&  
\qquad + \Delta t_n \left((I^h w^n)\,(I^h-I)\psi^n + \psi^n\,(I^h-I)w^n,Z^n_\rho \right)
\nonumber \\
& \leq   C\,\Delta t_n \left[ 
|I^h\psi^n-\Psi^n|_0\,|I^h w^n|_{0,\infty}\,|Z^n|_1 + h^2\,|I^h \psi^n|_1 \,|I^h w^n|_1 \, 
|Z^n|_{1,\infty} \right] 
\nonumber \\
&    \qquad + C\,\Delta t_n\,h \left[\, |w^n|_{0,\infty}\,|\psi^n|_1+
|\psi^n|_{0,\infty}\,|w^n|_1
\right]\,|Z^n|_1 
\nonumber \\
& \leq   C\,\Delta t_n
\left[  \, \left|D_t\vE^n
%\frac{\vE^n-\vE^{n-1}}{\Delta t_n}
\right|_0 + | \vE^n |_1 + h  %+ \Delta t_n   
\right] |Z^n|_1.  
\label{B4}
\end{align}
On noting (\ref{innerh}), (\ref{lumpev}), (\ref{invh},b), (\ref{hqu}),
Sobolev embedding, (\ref{vxwnreg}) and
(\ref{vpsierr}), we have that
\begin{align}
B_5 & = - \Delta t_n \left( I^h \psi^n \, Z^n, Z^n_\rho \right)^h - \Delta t_n 
\left( ( \Psi^n - I^h \psi^n) \, Z^n, Z^n_\rho \right)^h  \nonumber \\
& \leq C \,\Delta t_n \, 
|\psi^n|_{0,\infty}\,
| Z^n |_0 \, | Z^n |_1 + C\, \Delta t_n\,  | \Psi^n - I^h \psi^n |_0 \, 
| Z^n |_{0,\infty} \, | Z^n |_1  \nonumber \\
& \leq C \,\Delta t_n\, | Z^n |_0 \, | Z^n |_1 + C \,\Delta t_n \left[ 
| \vE^n |_1 + \left|D_t\vE^n
%\frac{\vE^n-\vE^{n-1}}{\Delta t_n}  
\right|_0 + h \right]
h^{-\frac{1}{2}} \,| Z^n |_0 \, | Z^n |_1 \nonumber \\
& \leq C \,\Delta t_n\,  |Z^n |_0 \, | Z^n |_1
\left[ 1 + h^{-\frac{1}{2}}\left(| \vE^n |_1 + 
\left| D_t\vE^n
%\frac{\vE^n-\vE^{n-1}}{\Delta t_n} 
\right|_0 \right) \right]. 
%+ C\, h^{-\frac{1}{2}} \,| Z^n |_0 \, 
%\Delta t_n \, \left|\frac{\vE^n-\vE^{n-1}}{\Delta t_n}\right|_0  \,
%| Z^n |_1.  
\label{B5}
\end{align}
Finally, we deduce from  (\ref{lumpev},b) and (\ref{Ih}) that  
\begin{eqnarray}
\lefteqn{
B_6  =
\Delta t_n\left[ \left(\left[|\vx_\rho^n|-|\vX^n_\rho|\right]\,g(v^n,w^n),Z^n\right)
+ \left(|\vX^n_\rho|\,(I-I^h)g(v^n,w^n),Z^n\right) \right] }
\nonumber \\
& &  + \Delta t_n 
\left[
\left( |\vX_\rho^n|\,I^h\left[g(v^n,w^{n})\right],Z^n\right)-
\left( |\vX_\rho^n|\,I^h\left[g(v^n,w^{n})\right],Z^n\right)^h\right]
\nonumber \\
& & + \Delta t_n 
\left[\left( |\vX_\rho^n|\left[g(I^h v^n,I^h w^{n})-g(V^n,W^{n-1})\right],Z^n\right)^h\right]
\nonumber \\
& \qquad \leq &  C\,\Delta t_n \left[ |\vx^n-\vX^n|_1\,|g(v^n,w^n)|_{0,\infty}
%\Delta t_n\left[ \left(|\vx_\rho^n|\,g(v^n,w^n),Z^n\right)
%-\left( |\vX_\rho^n|\,g(V^n,W^{n-1}),Z^n\right)^h\right]
+ |(I-I^h)g(v^n,w^n)|_0\right]|Z^n|_0
\nonumber \\
& & + C\,\Delta t_n \left[ h \,  |I^h[g(v^n,w^{n})]|_1 
+ \sqrt{(r^n,r^n)^h} \right] | Z^n |_0, \label{B6a}
\end{eqnarray}
where we have set $r^n = g(I^h v^n,I^h w^n)-g(V^n,W^{n-1})$. Using (\ref{fgreg},b) together 
with the fact that $| W^{n-1} |_{0,\infty} \leq C$, recall (\ref{disupw}),
we infer that
\begin{align*}
 | r^n | & \leq  | g(I^h v^n,I^h w^n) - g(I^h v^n,W^{n-1})| + | g(I^h v^n,W^{n-1}) 
 - g(V^n,W^{n-1})|  \\
& \leq C \left[ | I^h w^n - W^{n-1} | + | I^h v^n - V^n| \right],
\end{align*}
from which we deduce that $\sqrt{(r^n,r^n)^h} \leq C \bigl[\, | I^h w^n - W^{n-1} |_0 
+ | I^h v^n - V^n |_0 \bigr]$
on noting (\ref{innerh}) and (\ref{lumpev}).
Inserting this bound into (\ref{B6a}) and recalling (\ref{vxwnreg}), (\ref{wreg}), 
(\ref{Ih}),
(\ref{vpsierr}) 
as well as $\Delta t_n \leq Ch$ we have
\begin{equation}  \label{B6}
B_6 \leq 
  C\,\Delta t_n \left[\,  %\,\left(\|v^n\|_1 + \|w^n\|_1 + \|w^{n-1}\|_1
%\right) 
%+ \Delta t_n %\, \|w_t\|_{L^\infty(t_{n-1},t_n;L^2(\Ip))}
%\nonumber \\
%& \hspace{1in}
|\vE^n|_1 + |Z^{n-1}|_0 + \left|D_t\vE^n
%\frac{\vE^n-\vE^{n-1}}{\Delta t_n}
\right|_0 +h \right]
|Z^n|_0.
\end{equation}

It remains to analyse $B_2$. We claim that:
\begin{align}
|B_2| & \leq  C\, \Delta t_n  
\left[ \,\left|%\frac{\vE^n-\vE^{n-1}}{\Delta t_n}
D_t\vE^n \right|_0 + | \vE^n|_1 
+ | \vE^{n-1} |_1 +  h %+ \Delta t_n 
\right] \| Z^n \|_1   \nonumber  \\
& \qquad \qquad + C \, \Delta t_n  
\left|D_t\vE^n
%\frac{\vE^n-\vE^{n-1}}{\Delta t_n}
\right|^2_0 \,h^{-\frac{1}{2}}\, | Z^n |_0
+  C \,| Z^n |_{0,\infty} \, | \vE^n - \vE^{n-1} |_1^2 \nonumber \\
& \qquad \qquad 
+ C\,h \left( \Delta t_n \int_{t_{n-1}}^{t_n} \zeta(t) dt \right)^{\frac{1}{2}}  \Vert Z^n \Vert_1.
\label{B2}
\end{align}
In order to prove (\ref{B2}), we first write
\begin{eqnarray}
\lefteqn{ \hspace{-1cm}
B_2 = \left[\left( |\vX_\rho^n|\, I^h w^n-|\vX_\rho^{n-1}|\,I^h w^{n-1},Z^n \right)^h-
\left( |\vX_\rho^n|\, I^h w^n-|\vX_\rho^{n-1}|\,I^h w^{n-1},Z^n \right) \right] } 
\nonumber \\
& & \qquad   +
\left( |\vX_\rho^n|\, (I^h-I)\,w^n-|\vX_\rho^{n-1}|\,(I^h-I)\,w^{n-1},Z^n \right)
\nonumber \\
& & \qquad +
\left[\left( |\vx_\rho^n|\,w^n-|\vx_\rho^{n-1}|\,w^{n-1},Z^n \right)
- \Delta t_n \left( \left(|\vx_\rho|\,w\right)_t(t_n),Z^n\right)\right] \nonumber \\
& & \qquad +
\left( \left[|\vX_\rho^n|-|\vx_\rho^n|\right]\,w^n-\left[|\vX_\rho^{n-1}|
-|\vx_\rho^{n-1}|\right]\,w^{n-1},Z^n \right)
 =: \sum_{i=1}^4 B_{2,i}.
\label{B2a}
\end{eqnarray}
We now bound the $B _{2,i}$ terms defined in (\ref{B2a}).
It follows from  (\ref{lumperr}), (\ref {disnodegen}), (\ref{Ih}), (\ref{vxwnreg}), (\ref{DtX})
and as $w_t \in L^{\infty}(0,T;H^1(\Ip))$
%Not worth it, as we need it elsewhere. \label{hqua}} 
%and (\ref{disnodegen}) 
that
\begin{align}
|B_{2,1}| &\leq C\,h\,|\,|\vX_\rho^n|\, I^h w^n-|\vX_\rho^{n-1}|\,I^h w^{n-1}|_0\,|Z^n|_1 \nonumber \\
& \leq C\,h \left[|\vX^n-\vX^{n-1}|_1\,|w^n|_{0,\infty} +|I^h (w^n-w^{n-1})|_0\right]\,|Z^n|_1
\nonumber \\
&\leq C\,\Delta t_n\left[\,\left|D_t\vE^n\right|_0 + h 
\right]|Z^n|_1.
\label{B11}
\end{align} 
Similarly, we deduce from (\ref {disnodegen}), Sobolev embedding, 
(\ref{Ih}), (\ref{wreg}) and (\ref{DtX})  
%(\ref{hqu}) %\footnote{Same as footnote \ref{hqua}.} 
that
\begin{align}
|B_{2,2}| 
&\leq  C\left[ |\vX^n -\vX^{n-1}|_1 \,|(I-I^h)w^n|_{0,\infty}
+ |(I-I^h)(w^{n}-w^{n-1})|_0
\right]
|Z^n|_0 
\nonumber \\
&\leq  C\,h\left[ |\vX^n -\vX^{n-1}|_1 \,|w^n|_2
+ |w^{n}-w^{n-1}|_1
\right]
|Z^n|_0 
\nonumber \\
&\leq C \,\Delta t_n \left[\, \left|%\frac{\vE^n-\vE^{n-1}}{\Delta t_n}
D_t\vE^n \right|_0  + h  
\right]|Z^n|_0.
\label{B12}
\end{align}
In view of (\ref{vxreg},d), we have that 
$\vx_t \in L^{\infty}(0,T;[W^{1,\infty}(\Ip)]^2)$, $w_t \in L^{\infty}(0,T;L^{\infty}(\Ip))$
and hence, as $\Delta t_n \leq C\,h$,
\begin{align}
|B_{2,3}| &\leq C\,\Delta t_n \left( \Delta t_n \int_{t_{n-1}}^{t_n}|(|\vx_\rho|\,w)_{tt}|^2_0 \dt \right)^{\frac{1}{2}}
|Z^n|_0
\nonumber \\
&\leq C\,h \left( \Delta t_n \int_{t_{n-1}}^{t_n}
\left[ 1+
|w_{tt}|^2_0  + |\vx_{tt}|_{1}^2
\right]\dt \right)^{\frac{1}{2}}
|Z^n|_0.
\label{B14}
\end{align}
It remains to analyse $B_{2,4}$, which requires some more intricate arguments. To begin,
\begin{eqnarray}
B_{2,4} &= &
\left( \left[\vX_\rho^n \cdot \vTau^n -\vx_\rho^n \cdot \vtau^n \right]\,w^n-
\left[\vX_\rho^{n-1} \cdot \vTau^{n-1} -\vx_\rho^{n-1} \cdot \vtau^{n-1} \right]\,w^{n-1},Z^n \right)
\nonumber \\
&= & \left(\vX_\rho^n \cdot (\vTau^n-\vtau^n)\,w^n
-\vX_\rho^{n-1} \cdot (\vTau^{n-1}-\vtau^{n-1})\,w^{n-1},Z^n\right)
\nonumber \\
& &  \qquad +\left((\vX^n-\vx^n)_\rho \cdot \vtau^n\, w^n -(\vX^{n-1}-\vx^{n-1})_\rho 
\cdot \vtau^{n-1} \,w^{n-1},Z^n \right) \nonumber \\
&  =: &  B_{2,4,1}+ B_{2,4,2}. \label{B13}
\end{eqnarray}
On recalling (\ref{disctn}), we rewrite
\begin{eqnarray}
B_{2,4,1} & = & \tfrac{1}{2}\left(|\vX_\rho^n|\,|\vtau^n-\vTau^n|^2\,w^n
-|\vX_\rho^{n-1}|\,|\vtau^{n-1}-\vTau^{n-1}|^2\,w^{n-1},Z^n\right) 
\nonumber \\
& = & \tfrac{1}{2}\left( \left[ |\vX_\rho^n|- |\vX_\rho^{n-1} |\right] \,|\vtau^n-\vTau^n|^2\,w^n, 
Z^n\right)   \nonumber \\
& & 
+ \tfrac{1}{2}\left(|\vX_\rho^{n-1}|\, |\vtau^{n-1} -\vTau^{n-1}|^2 \,(w^n - w^{n-1}), Z^n 
\right) \nonumber \\
&  &   
+ \tfrac{1}{2}\left(|\vX_\rho^{n-1} |\left[ |\vtau^n-\vTau^n|^2 - |\vtau^{n-1} 
-\vTau^{n-1}|^2 \right] w^n, Z^n \right)  
=: I + II + III.
\label{B131}
\end{eqnarray}
%Using (\ref{invh}), 
Noting (\ref{DtX}),
(\ref{tauerr}), 
 Sobolev embedding  
 %$H^1(\Ip) \hookrightarrow C^0(\Ip)$, 
%(\ref{invh}), (\ref{hqu}) 
and (\ref{Enh}), we obtain that
\begin{align}
| I | & \leq C\, | \vX^n - \vX^{n-1} |_1 \, | \vtau^n - \vTau^n |_0 \, 
| \vtau^n - \vTau^n |_{0,\infty} \, | Z^n |_{0,\infty} \nonumber \\
& \leq C \Delta t_n  \left[ h^{-1}\, | D_t \vE^n |_0 + 1 \right] 
\,  h^{-\frac{1}{2}} \,\left[| \vE^n |_1 + h  \right]^2 \Vert  Z^n \Vert_1  \nonumber \\
& \leq C \,\Delta t_n \, \left[  \left|
D_t\vE^n\right|_0 + h \right] \Vert Z^n \Vert_1. \label{termI} 
\end{align}
Since $w_t \in L^{\infty}(0,T;L^{\infty}(\Ip))$ and $| \vtau^{n-1} - \vTau^{n-1} |_{0,\infty} 
\leq 2$, we deduce from (\ref{disnodegen}) and (\ref{tauerrnm})   
that 
%, and %similarly to (\ref{tauerr}) noting (\ref{disnodegen})
%\footnote{We should introduce the $(n-1)$ version of (\ref{tauerr}) 
%after (\ref{disnodegen}) and refer to that. This can also be used in (\ref{A134}).
%\label{tauerr-1}}, 
%we infer that
\begin{equation} 
\label{termII} \qquad \qquad | II |  \leq C \,\Delta t_n \,| \vtau^{n-1} - \vTau^{n-1} |_0 \, 
| Z^n |_0 \leq C\, \Delta t_n \left[\, | \vE^{n-1} |_1 + h \right] |Z^n |_0.
\end{equation}

In order to treat term $III$, we first observe that
\begin{align}
|\vtau^n-\vTau^n|^2 - |\vtau^{n-1} -\vTau^{n-1}|^2 &= 
- (\vtau^n - \vtau^{n-1}) \cdot ( \vTau^n + \vTau^{n-1} ) \nonumber \\
& \qquad -( \vTau^n - \vTau^{n-1}) \cdot ( \vtau^{n}+\vtau^{n-1}). \label{termIIIc}
\end{align}
Furthermore, a straightforward calculation shows that
\begin{subequations}
\begin{align}
\vTau^n - \vTau^{n-1} & = \frac{1}{| \vX^{n-1}_\rho |} 
\left[  I - \vTau^n \otimes \frac{\vX^{n}_\rho + \vX^{n-1}_\rho}{
| \vX^{n}_\rho | + | \vX^{n-1}_\rho | } \right] ( \vX^n_\rho - \vX^{n-1}_\rho)  \nonumber \\
& =   \frac{1}{| \vX^{n-1}_\rho |} (P^n+R^n)\,( \vX^n_\rho - \vX^{n-1}_\rho), 
\label{termIIIa}
\end{align}
 where
 \begin{align}
 P^n= I - \vTau^n \otimes \vTau^n \quad \mbox{ and } \quad
| R^n | \leq C \,| \vX^n_\rho - \vX^{n-1}_\rho |.
\label{termIIIb}
\end{align}
\end{subequations}
Inserting (\ref{termIIIc}) and (\ref{termIIIa},b) into $III$ and observing that 
$(\vtau^n - \vtau^{n-1}) \cdot (\vtau^n + \vtau^{n-1})=0$, we derive
\begin{align}
 2 \, III   &=  - \left( | \vX^{n-1}_\rho |\, (\vtau^n - \vtau^{n-1}) 
 \cdot (\vTau^n  + \vTau^{n-1})\, w^n, Z^n \right)  \nonumber  \\
 & \qquad  - \left( (P^n+R^n) \,(\vX^n_\rho - \vX^{n-1}_\rho) 
 \cdot (\vtau^{n}+\vtau^{n-1}) \, w^n,Z^n \right) \nonumber \\
  & = - \left( | \vX^{n-1}_\rho | \,(\vtau^n - \vtau^{n-1}) 
  \cdot (\vTau^n - \vtau^n + \vTau^{n-1} -\vtau^{n-1})\, w^n, Z^n \right)  
  \nonumber \\
 & \qquad    - \left( \left[ (I^h \vx^n)_\rho - (I^h \vx^{n-1})_\rho \right] 
 \cdot \left( P^n+(R^n)^t \right) (\vtau^{n} + \vtau^{n-1})\, w^n, Z^n \right) 
 \nonumber \\
& \qquad   + \left( P^n ( \vE^n_\rho - \vE^{n-1}_\rho) \cdot (\vtau^{n} + \vtau^{n-1})\, w^n, 
Z^n \right)  \nonumber 
\\
& \qquad +  \left( R^n ( \vE^n_\rho - \vE^{n-1}_\rho) \cdot (\vtau^{n} + \vtau^{n-1})\, w^n, 
Z^n \right) \nonumber \\  
&=:   III_1 + III_2 + III_3 + III_4. 
\label{termIII} 
\end{align}
It follows from (\ref{disnodegen}), (\ref{tauerrnm}), (\ref{tauerr}), Sobolev embedding and 
(\ref{vxreg},d)
%\footnote{The assumed 
%$\vx \in W^{1,\infty}(0,T;[H^2(\Ip)]^2)$ in (\ref{vxreg}) yields that 
%$|\vtau^n-\vtau^{n-1}|_{0,\infty} 
%\leq C\,\Delta t_n$. State this.} 
that
\begin{eqnarray}
| III_1 | &\leq &  |\vtau^n-\vtau^{n-1}|_{0,\infty}\left[ |\vtau^n-\vTau^n|_0 + |\vtau^{n-1}-\vTau^{n-1}|_0
\right] |w^n|_{0,\infty}\, |Z^n|_0
\nonumber \\
&\leq &  C \,\Delta t_n %\|\vx_t\|_{L^\infty(t_{n-1},t_n; [W^{1,\infty}(\Ip)]^2)}
\left[\, | \vE^n |_1 + | \vE^{n-1} |_1 +h \right] | Z^n |_0,
\label{termIII_1}
\end{eqnarray}
as $\vx_t \in L^\infty(t_{n-1},t_n; [W^{1,\infty}(\Ip)]^2)$.
Recalling the definition of $P^n$ we may write
\begin{align*}
& P^n (\vtau^{n} + \vtau^{n-1})  =  (\vtau^{n} - \vTau^n)  
+ (\vtau^{n-1}  -  \vTau^n)  + \tfrac{1}{2}\left[  | \vtau^{n} - \vTau^n |^2 + 
|\vtau^{n-1} - \vTau^n |^2 \right] \vTau^n,
\end{align*}
so, similarly to (\ref{termIII_1}), we have from (\ref{termIIIb}), (\ref{tauerrnm}),
(\ref{tauerr}), Sobolev embedding and (\ref{vxreg},d) that
%\footnote{Need to add some other references here.}  
\begin{align}
 | III_2 |  & \leq C\, \Delta t_n \, \left[
 |\vtau^n-\vTau^n|_0 + |\vtau^{n-1}-\vTau^{n-1}|_0 + |\vX^n-\vX^{n-1}|_1
 \right]\,|w^n|_{0,\infty}\, |Z^n|_0 
\nonumber \\
 & \leq C  \,\Delta t_n  \left[\, | \vE^n |_1 + 
| \vE^{n-1} |_1 +h \right] | Z^n |_0,  %+ C \Delta t_n 
%\left|\frac{\vE^n-\vE^{n-1}}{\Delta t_n}\right|_0 | Z^n |_0.
%\nonumber
\label{termIII_2}
\end{align}
since $\Delta t_n \leq C\,h$.
Next, performing integration by parts on the subintervals $\sigma_j$ we derive that
\begin{align}
III_3  & = - \left( P^n ( \vE^n - \vE^{n-1}), %\cdot  
\left[ (\vtau^{n} + \vtau^{n-1})\, w^n \right]_\rho \, Z^n %\right) \nonumber \\
%& \qquad - \left( P^n ( \vE^n - \vE^{n-1}) \cdot 
+ (\vtau^{n-1} + \vtau^n)\, w^n\, Z^n_\rho \right) 
\nonumber \\
& \qquad - \sum_{j=1}^J \left( P^n \mid_{\sigma_{j+1}} - P^n\mid_{\sigma_j} \right) 
( \vE^n - \vE^{n-1})(\rho_j) \cdot
( \vtau^{n}+ \vtau^{n-1})\, w^n(\rho_j) \,Z^n(\rho_j),
\label{termIII_3a}
\end{align}
where $\sigma_{J+1} \equiv \sigma_1$.
Introducing the nodal basis functions $\chi_j \in V^h_1$ such that $\chi_j(\rho_i) = \delta_{i,j}$,
$i,\,j=1,\ldots,J$, and then choosing $\vxi^h=\chi_j\,\vxi$, for any $\vxi \in \bR^2$, in 
%Inserting the nodal basis functions $\chi_j,j=1,\ldots,J$ into 
(\ref{G1wbn}) we obtain 
\begin{align}
&  \left[\vX^n_{\rho} \mid_{\sigma_{j+1}} - \vX^n_{\rho} \mid_{\sigma_j}\right] \cdot \vxi =  
-  \left( |\vX^{n-1}_\rho|^2\,f(W^{n-1})  
\,\vNu^{n-1} , \chi_j\,\vxi\right)^h  \nonumber \\
& \hspace{1in} +  %\frac{1}{\Delta t_n} 
\left( |\vX^{n-1}_\rho|^2 \left[\alpha \, D_t \vX^n %(\vX^n-\vX^{n-1}) 
+ (1-\alpha) 
\left(D_t \vX^n
%(\vX^n-\vX^{n-1}) 
\cdot \vNu^{n-1} \right) 
\vNu^{n-1} \right] , \chi_j \,\vxi\right)^h. \nonumber
%\label{G1wbnj}
\end{align}
We deduce from this and (\ref{termIIIb}), with the help of (\ref{disnodegen}), (\ref{disupw}),
Sobolev embedding 
and (\ref{vxreg}), 
that for $j=1,\ldots,J$
%\footnote{We require $\vx \in W^{1,\infty}(0,T;[L^\infty(\Ip)]^2)$ here, 
%and this follows from the assumed $\vx \in W^{1,\infty}(0,T;[H^2(\Ip)]^2)$ in (\ref{vxreg}). State this.}
\begin{align}
\left| P^n \mid_{\sigma_{j+1}}\! - P^n \mid_{ \sigma_j} \right|  &\leq C \left|  \vX^n_{\rho}  
 \mid_{\sigma_{j+1}}\! - \vX^n_{\rho} \mid_{\sigma_j} \right|    \nonumber \\
&\leq   C \,h \left[   | f(W^{n-1}) |_{0,\infty} + \left|  (D_t\vX^n)(\rho_j)
 %\frac{ \vX^n(\rho_j) - \vX^{n-1}(\rho_j)}{\Delta t_n} 
 \right|  \right]   \nonumber \\
 & \leq  C \,h \left[ 1 +  
 \left| %\frac{\vx^n-\vx^{n-1}}{\Delta t_n}
 D_t\vx^n
 \right|_{C(\Ip)} +
 \left|  %\frac{ \vE^n(\rho_j) - \vE^{n-1}(\rho_j)}{\Delta t_n} 
 (D_t\vE^n)(\rho_j)\right|
 \right]\nonumber \\ 
 &\leq C \,h \left[ 1 +   
 \left|(D_t\vE^n)(\rho_j)\right|
 \right]. \nonumber
 %\label{termIII_3b}
 \end{align}
Inserting this estimate into (\ref{termIII_3a}) and recalling 
(\ref{termIIIb}), (\ref{nodegen}),
(\ref{hqu}), (\ref{lumpev}), Sobolev embedding, (\ref{vxwnreg})
and (\ref{invh}), we deduce that
\begin{align}
| III_3 | & \leq C \,\Delta t_n 
\left|
D_t\vE^n
%\frac{\vE^n-\vE^{n-1}}{\Delta t_n}
\right|_0 \,\|Z^n \|_1 + C\, \Delta t_n 
\left|
D_t\vE^n
%\frac{\vE^n-\vE^{n-1}}{\Delta t_n}
\right|^2_0 \,| Z^n |_{0,\infty}.
\nonumber \\
 & \leq C \,\Delta t_n 
\left|
D_t\vE^n
\right|_0 \,\|Z^n \|_1 + C\, \Delta t_n 
\left|
D_t\vE^n
%\frac{\vE^n-\vE^{n-1}}{\Delta t_n}
\right|^2_0 \,h^{-\frac{1}{2}}\, | Z^n |_0.
\label{termIII_3}
\end{align}
Finally, we have, on noting (\ref{termIIIb}), (\ref{A11}), Sobolev embedding
and (\ref{vxwnreg}), that
\begin{align}
| III_4 | & \leq C \,| \vE^n -\vE^{n-1} |_1 \left[  
| \vE^n -\vE^{n-1} |_1 + | I^h(  \vx^n - \vx^{n-1}) |_1 \right] \,|w^n|_{0,\infty} 
%| \vE^n -\vE^{n-1}  |_1
\, | Z^n |_{0,\infty} \nonumber \\
& \leq  C \,| Z^n |_{0,\infty}  \,| \vE^n - \vE^{n-1} |_1^2  +
C \,\Delta t_n  \, | \vE^n - \vE^{n-1} |_1 \,  \|Z^n \|_1.
\label{termIII_4}
\end{align}
Combining (\ref{B131})--( \ref{termII}), 
(\ref{termIII})--(\ref{termIII_2}), (\ref{termIII_3}) and (\ref{termIII_4}),
we conclude that 
\begin{align}
| B_{2,4,1} | & \leq  C\, \Delta t_n  
\left[ \,\left|D_t\vE^n
%\frac{\vE^n-\vE^{n-1}}{\Delta t_n}
\right|_0 + | \vE^n |_1 
+ | \vE^{n-1} |_1 +  h %+ \Delta t_n 
\right] \| Z^n \|_1    \nonumber\\
& \quad + C \, \Delta t_n  
\left|
D_t\vE^n
%\frac{\vE^n-\vE^{n-1}}{\Delta t_n}
\right|^2_0 \,h^{-\frac{1}{2}} \,| Z^n |_0
+  C \,| Z^n |_{0,\infty} \, | \vE^n - \vE^{n-1} |_1^2.
\label{b131est}
\end{align}
For the second term in (\ref{B13}) we obtain with the help of integration by parts that
\begin{align*}
&B_{2,4,2} = \Delta t_n  \left( D_t[(\vX^n-\vx^n)_{\rho} \cdot \tau^n \, w^n],Z^n \right)  \\
& \quad =   \Delta t_n \left(D_t [(\vX^n - \vx^n)_{\rho}] \cdot \tau^n w^n, Z^n \right) 
+ \Delta t_n \left( (\vX^{n-1} - \vx^{n-1})_{\rho} \cdot D_t [ \tau^n w^n], Z^n \right) \\
& \quad =  - \Delta t_n \left( D_t (\vX^n -\vx^n), (w^n\,Z^n\,\vtau^n)_\rho\right) 
+ \Delta t_n \left( (\vX^{n-1} - \vx^{n-1})_{\rho} \cdot D_t [ \tau^n w^n], Z^n \right).
\end{align*}
Similarly to (\ref{termIII_1}), it follows from  
(\ref{Ih}), (\ref{vxreg}--e) and Sobolev embedding that
\begin{align}
|B_{2,4,2}| & \leq   C\, \Delta t_n \,\left[ |(I-I^h) D_t \vx^n |_0 + \Delta t_n \left|D_t\vE^n\right|_0\right]
\, \|Z^n\|_1  \nonumber \\
&   \qquad \qquad + C\, \Delta t_n \, \left[\, |(I-I^h)\vx^{n-1}|_1 + |E^{n-1}|_1 \right]  \,|Z^n|_0
\nonumber \\
& \leq    C\,
\Delta t_n\left[ \,
\left|D_t\vE^n
%\frac{\vE^n-\vE^{n-1}}{\Delta t_n}
\right|_0  +  |\vE^{n-1}|_1 +h  \right] \|Z^n\|_1.
\label{B132}
\end{align}
\begin{comment}
\left|\,\left[|\vX_\rho^n|-|\vx_\rho^n|\right]\,w^n-\left[|\vX_\rho^{n-1}|
-|\vx_\rho^{n-1}|\right]\,w^{n-1}\right|_0 \,|Z^n|_0
\nonumber \\
&\leq C\left[
\left|\,\left[|\vX_\rho^n|-|\vx_\rho^n|\right]
- \left[|\vX_\rho^{n-1}|-|\vx_\rho^{n-1}|\right]
\right|_0
+ |\vx^{n-1}-\vX^{n-1}|_1\,|w^n-w^{n-1}|_{0,\infty}
\right]
|Z^n|_0
\end{comment}
Thus, (\ref{B13}), (\ref{b131est}) and (\ref{B132}) yield the bound on $B_{2,4}$, 
which together with 
%(\ref{B11}), (\ref{B12}), (\ref{B14}) and 
(\ref{B2a})--(\ref{B14}) imply (\ref{B2}) on recalling (\ref{zeta}).

Hence we have bounded all the terms on the right hand side of (\ref{Zn}), so we obtain from
(\ref{B1})--(\ref{B5}), (\ref{B6}) and (\ref{B2})  
%on noting  (\ref{disnodegenk}), (\ref{invh}), (\ref{zeta}) 
%as well as $\Delta t_n \leq C\,h$, 
that
\begin{align}
%\lefteqn{
&\tfrac{1}{2}\,(|\vX^n_\rho|\,Z^n,Z^n)^h +  \frac{d}{2M} \, 
\Delta t_n  \, | Z^n |_1^2  + \tfrac{1}{2} \,(|\vX^{n-1}_\rho|\,(Z^n-Z^{n-1}),(Z^n-Z^{n-1})\,)^h 
%}
\nonumber \\
& \quad\leq   \tfrac{1}{2} \,(|\vX^{n-1}_\rho|\,Z^{n-1},Z^{n-1})^h +
\delta  \, \Delta t_n\,|Z^n |_1^2 + C(\delta)\,\Delta t_n
\left[\, | Z^n |_0^2 + |Z^{n-1}|_0^2\,  \right] \nonumber \\
&  \qquad \quad   + C(\delta) \, h^2 \int_{t_{n-1}}^{t_n} \zeta(t) \dt 
+ C(\delta)\, \Delta t_n \left[ |\vE^n|_1^2 + | \vE^{n-1} |_1^2 +  
\left|D_t\vE^n \right|^2_0 \right] \nonumber  \\
&  \qquad  \quad + C(\delta) \,\Delta t_n\,  
h^{-1} \,| Z^n |_0^2 \, \left[ |\vE^n|_1^2 +  
\left|D_t\vE^n \right|^2_0 \right] + C \,| Z^n |_{0,\infty}\, | \vE^n - \vE^{n-1} |_1^2. 
%\nonumber\\
\label{Zn1}
\end{align} 
Let us examine the last two terms appearing on the right hand side of 
(\ref{Zn1}).  We deduce from (\ref{Z0}), (\ref{Enh}), (\ref{invh}) and (\ref{hqu}) that 
\begin{subequations}
\begin{align}
& \Delta t_n \, h^{-1} \,| Z^n |_0^2 \left[ |\vE^n|_1^2 + \left| D_t\vE^n
\right|^2_0 \right] \nonumber \\
& \hspace{1in} \leq %\bigl( | Z^{n-1} |_{0,\infty} + C h^{-\frac{1}{2}} | Z^n - Z^{n-1} |_0  + 2 
C\,h^{-1} \left[ | Z^{n-1}  |_0^2 +  %2 h^{-1} 
| Z^n - Z^{n-1}|_0^2 \right] \Delta t_n \left[ |\vE^n|_1^2  +
\left|
D_t\vE^n
%\frac{\vE^n-\vE^{n-1}}{\Delta t_n}
\right|^2_0 \right] \nonumber \\
& \hspace{1in} \leq 
%\frac{C}{\beta^2} \,h \,e^{\gamma K}\,
C\,\Delta t_n 
\left[ |\vE^n |^2_1+ \left|D_t\vE^n
%\frac{\vE^n-\vE^{n-1}}{\Delta t_n}
\right|^2_0 \right]
%\nonumber \\
%& \qquad \qquad +  C \,%h \,e^{\gamma K}
+C\, h^{\frac{1}{2}} \, | Z^n- Z^{n-1} |_0^2, 
\label{Zn1e} 
\end{align}
and %finally with the help of (\ref{En4})
\begin{align}
%\lefteqn{ \hspace{-1cm}
| Z^n |_{0,\infty} \,| \vE^n - \vE^{n-1} |_1^2 &\leq C \,h^{-\frac{1}{2}}
\left[ | Z^{n-1} |_{0} 
+ %h^{-\frac{1}{2}}\, 
| Z^n - Z^{n-1} |_0 \right] | \vE^n - \vE^{n-1} |_1^2   \nonumber \\
& \leq C \, 
| \vE^n - \vE^{n-1} |_1^2 + C\, h^{\frac{1}{4}}\, | Z^n - Z^{n-1} |_0 
| \vE^n - \vE^{n-1} |_1  \nonumber \\
&\leq  C \, | \vE^n - \vE^{n-1} |_1^2
+ C \,h^{\frac{1}{2}}\,  | Z^n - Z^{n-1} |_0^2.
\label{Zn1f}
\end{align}
\end{subequations}
Combining (\ref{Zn1}) and (\ref{Zn1e},b)
  and choosing $\delta$ sufficiently small
yields, on noting (\ref{disnodegen}) and  (\ref{lumpev}),
\begin{eqnarray}
\lefteqn{ \hspace{-1cm}
 (|\vX^n_\rho|\,Z^n,Z^n)^h   + \frac{d}{2M} \, \Delta t_n  \, | Z^n |_1^2 +  
(1- C\,h^{\frac{1}{2}}) \, (|\vX^{n-1}_\rho|\,(Z^n-Z^{n-1}),(Z^n-Z^{n-1})\,)^h } \nonumber \\
& \leq & (|\vX^{n-1}_\rho|\,Z^{n-1},Z^{n-1})^h +  C\, \Delta t_n  
\left[ (|\vX^{n-1}_\rho|\,Z^{n-1},Z^{n-1})^h  + | \vE^{n-1} |_1^2 
%+  h^2 %+ (\Delta t_n)^2 
\right] \nonumber \\
& &\qquad \qquad + C \,\Delta t_n \left|D_t \vE^n %\frac{\vE^n-\vE^{n-1}}{\Delta t_n}
\right|^2_0 
+ C \,| \vE^n - \vE^{n-1} |_1^2 +  C \, h^2 \int_{t_{n-1}}^{t_n} \zeta(t) \dt.
\label{Zn2}
\end{eqnarray}
We now proceed by choosing $h^{\star}>0$ so small that $C\, (h^\star)^{\frac{1}{2}}\leq 1$ in (\ref{Zn2}).
Multiplying (\ref{Zn2}) by $\beta^2$ and adding to (\ref{En1})
%+ $\beta^2 \times$ (\ref{Zn2}) we obtain
yields that
\begin{eqnarray}
\lefteqn{
\Delta t_n \,\alpha \, \frac{m^2}{4}  |D_t\vE^n|_0^2
%\frac{\vE^n-\vE^{n-1}}{\Delta t_n},
%\frac{\vE^n-\vE^{n-1}}{\Delta t_n}
+ |\vE^n-\vE^{n-1}|_1^2
+ |\vE^n|_1^2 + \beta^2\, (|\vX^n_\rho|\,Z^n,Z^n)^h }  \nonumber\\
& \leq & | \vE^{n-1}|_1^2 + \beta^2 \,(|\vX^{n-1}_\rho|\,Z^{n-1},Z^{n-1})^h +
C \,\Delta t_n  \left[ | \vE^{n-1}|_1^2 +  (| \vX_\rho^{n-1} |\, Z^{n-1}, Z^{n-1})^h 
\right] \nonumber \\
& & \qquad  + C\,  h^2 \,  \int_{t_{n-1}}^{t_n} \zeta(t) \dt + C \,\beta^2 \,\Delta t_n 
\left|D_t\vE^n
%\frac{\vE^n-\vE^{n-1}}{\Delta t_n}
\right|^2_0 + C\, \beta^2\, | \vE^n - \vE^{n-1} |_1^2. 
\label{Zn2a}
\end{eqnarray}
Choosing $\beta \in (0,1]$ in such a way that $C\, \beta^2 \leq 
\min\{\alpha \frac{m^2}{4},1\}$ we find, 
with the help of the induction hypothesis (\ref{Errind2}), that for any $h \in (0,h^\star]$
\begin{align}
&|\vE^n|_1^2 + \beta^2 \,(|\vX^n_\rho|\,Z^n,Z^n)^h  
\nonumber \\
&\;\; \leq \left(1+\frac{C\,\Delta t_n}{\beta^2}\right)
\left[| \vE^{n-1}|_1^2 + \beta^2 \,(|\vX^{n-1}_\rho|\,Z^{n-1},Z^{n-1})^h \right] 
%& &\quad + \frac{C}{\beta^2}\, \Delta t_n  \left[ | \vE^{n-1}|_1^2 + 
%\beta^2 \,(| \vX_\rho^{n-1}\, | Z^{n-1}, Z^{n-1})^h \right] 
+ C\,  h^2 \int_{t_{n-1}}^{t_n} \zeta(t) \dt\nonumber \\
& \;\; \leq   h^2 \, e^{\gamma \int_0^{t_{n-1}} \zeta(t)\dt} 
\left[ 1 + C \left( 1 + \frac{1}{\beta^2} \right)
 \int_{t_{n-1}}^{t_n} \zeta(t) \dt \right] \leq h^2 \, e^{\gamma \int_0^{t_n} \zeta(t) \dt},
\label{Zn2b}
\end{align}
provided that $\gamma \geq C\left( 1 + \frac{1}{\beta^2} \right)$.
Hence, on assuming (\ref{Errind2}) for some $n \in \{1,\ldots,N\}$ we have shown that it holds for $n+1$,
(\ref{Zn2b}), provided that $\beta \in (0,1]$ is chosen as required in deriving (\ref{En3}) and (\ref{Zn2b}),
$\gamma$ as required by (\ref{gamma1}) and in deriving (\ref{Zn2b}), and finally $h^\star$     
so that (\ref{hstar}), (\ref{Psinbd})  and the condition for deriving (\ref{Zn2a}) are satisfied.

Therefore, under the above constraints on $\beta,\,\gamma$ and $h^\star$, 
we have that (\ref{Zn2b}) and (\ref{disnodegen}) hold for $n=0,\ldots,N$.
Moreover, we have that (\ref{En4}) and (\ref{Zn2}) hold for $n=1,\ldots,N$.
We deduce from (\ref{Zn2b}) and (\ref{disnodegen}), on noting (\ref{lumpev}), 
and on summing %the first inequality in 
(\ref{En4}) and (\ref{Zn2}) from $n=1$ to $N$ that
\begin{align}
\sup_{n=0,\ldots,N} \left[|\vE^n|_1^2 + |Z^n|_0^2\right] +
\sum_{n=1}^N \Delta t_n\, \left[ \left| D_t\vE^{n}\right|_0^2 + |Z^n|_1^2 \right]
\leq C\,h^2.
\label{errfinal}
\end{align} 
Finally, (\ref{errfinal}), (\ref{EnZn}), (\ref{Ih}) and (\ref{vxreg},d) yield the desired
error bounds (\ref{mainerr}) of Theorem \ref{main}.

\setcounter{equation}{0}
\section{Numerical results} \label{sec:5}
Throughout this section we set $d=1$.
Let us begin by investigating the experimental order of convergence (eoc). We use the following  example
from Section 6.4 of \cite{PozziS15} and consider for $\rho \in \Ip, \, t\in [0,1]$ 
$$
\vec{x}(\rho,t)=\left(\begin{array}{ll}(1+\frac12 \sin(2\pi t))\,\cos(2\pi \rho)\\[1mm]
(1-\frac12 \sin(2\pi t))\,\sin(2\pi \rho)\end{array}\right), \quad 
w(\rho,t)=t\,\cos(8\pi \rho) + (1-t)\,\sin(6\pi \rho),
$$ 
as well as $f(w)=2w$, $g(v,w)=0$. With these choices (\ref{xtang}), (\ref{tildew}) are satisfied
with additional terms $\vec{S}$, $S_w$ on the right hand side respectively, so that we
replace $f(W^{n-1}) \,\vNu^{n-1}$ by $f(W^{n-1})\,\vNu^{n-1} + \vec S$  in (\ref{G1wbn}) and
$g(V^{n},W^{n-1})$ by $g(V^{n},W^{n-1}) + S_w$ in (\ref{G2wbn}).
In what follows we choose a uniform 
element length $h=1/J$ and, unless stated otherwise, 
a uniform time step $\Delta t = h^2$ throughout  the computations in this section.  
We monitor the following errors:
%\footnote{Are these really calculated, or is 
%$w^n$ replaced by $I^h\,w^n$ etc. ?}
\begin{displaymath}
\mathcal{E}_1:=\sup_{n=0,\ldots,N} |Z^n|_0^2,\;
\mathcal{E}_2:=\sum_{n=1}^N \Delta t\, |Z^n|_1^2, \;
\mathcal{E}_3:=\sup_{n=0,\ldots,N} |\vE^n|_1^2, \;
\mathcal{E}_4:=\sum_{n=1}^{N}\Delta t\,
\left| D_t\vE^{n}\right|_0^2.
\end{displaymath}
In Tables \ref{table:1} and \ref{table:2} we display the values of $\mathcal{E}_i$, 
$i=1,\ldots,4$, for $\alpha = 1$ with $\Delta t = h^2$ and $\Delta t = 0.5h$ respectively. For 
$\Delta t = 0.5h$ we see eocs close to two, while for $\Delta t = h^2$ we see eocs close to four. 
When $\Delta t = h^2$ the eocs for $\mathcal{E}_2$ and $\mathcal{E}_3$ are better than the 
spatial approximation error, thus demonstrating superconvergence. Tables \ref{table:3} 
and \ref{table:4} show the corresponding results for $\alpha = 0.1$, again we see eocs 
close to two for $\Delta t = 0.5h$ and eocs close to four for $\Delta t = h^2$.
We saw similar results when we computed $\mathcal{E}_i$, $i=1,\ldots,4$, for the scheme 
in \cite{PozziS15}. The above results confirm the bounds obtained in Theorem \ref{main}.
It will be a subject of future research to rigorously prove the better rates in the case 
$\Delta t =h^2$.

 \begin{center}
\begin{table}[!h]
 \begin{tabular}{ |c||c|c|c|c|c|c|c|c| }
 \hline
J & $\mathcal{E}_1\times10$ & ${\rm eoc}_1$ & $\mathcal{E}_2$ & ${\rm eoc}_2$ & 
$\mathcal{E}_3$ & ${\rm eoc}_3$ & $\mathcal{E}_4\times10$ & ${\rm eoc}_4$ \\ 
 \hline
 \hline
 30 & 0.1817509 & - & 2.6729150 & - & 0.2043562 & - & 0.1177607 & -   \\ 
 60 & 0.0116485 & 3.96 & 0.1712537 & 3.96 & 0.0134930 & 3.92 & 0.0079150 & 3.90  \\ 
120 & 0.0007333 & 3.99 & 0.0107777 & 3.99 & 0.0008574 & 3.98 & 0.0005056 & 3.97  \\ 
240 & 0.0000459 & 4.00 & 0.0006748 & 4.00 & 0.0000538 & 3.99 & 0.0000318 & 3.99  \\ 
\hline
\end{tabular}
\caption{$\alpha = 1, \Delta t = h^2$}\label{table:1}
\end{table}
\end{center}
\begin{center}
\begin{table}[!h]
 \begin{tabular}{ |c||c|c|c|c|c|c|c|c| }
 \hline
J & $\mathcal{E}_1\times10$ & ${\rm eoc}_1$ & $\mathcal{E}_2$ & ${\rm eoc}_2$ & 
$\mathcal{E}_3$ & ${\rm eoc}_3$ & $\mathcal{E}_4\times10$ & ${\rm eoc}_4$ \\ 
 \hline
 \hline
 30 & 0.5093612 & - & 5.7947280 & - & 0.4464560 & - & 0.2664600 & -   \\ 
 60 & 0.0835780 & 2.61 & 0.9498929 & 2.61 & 0.0868606 & 2.36 & 0.0545729 & 2.29  \\ 
120 & 0.0162823 & 2.36 & 0.1859501 & 2.35 & 0.0187461 & 2.21 & 0.0125972 & 2.12  \\ 
240 & 0.0035578 & 2.19 & 0.0408490 & 2.19 & 0.0043234 & 2.12 & 0.0030531 & 2.04  \\ 
\hline
\end{tabular}
\caption{$\alpha = 1, \Delta t = 0.5h$}\label{table:2}
\end{table}
\end{center}
\begin{center}
\begin{table}[!h]
 \begin{tabular}{ |c||c|c|c|c|c|c|c|c| }
 \hline
J & $\mathcal{E}_1\times10$ & ${\rm eoc}_1$ & $\mathcal{E}_2$ & ${\rm eoc}_2$ & 
$\mathcal{E}_3$ & ${\rm eoc}_3$ & $\mathcal{E}_4\times10$ & ${\rm eoc}_4$ \\ 
 \hline
 \hline
 30 & 0.3227043 & - & 4.8893170 & - & 1.0377730 & - & 0.6239216 & -   \\ 
 60 & 0.0205566 & 3.97 & 0.2808019 & 4.12 & 0.0562308 & 4.21 & 0.0372623 & 4.07  \\ 
120 & 0.0012831 & 4.00 & 0.0172030 & 4.03 & 0.0033578 & 4.07 & 0.0022726 & 4.04  \\ 
240 & 0.0000801 & 4.00 & 0.0010698 & 4.01 & 0.0002074 & 4.02 & 0.0001411 & 4.01  \\ 
\hline
\end{tabular}
\caption{$\alpha = 0.1, \Delta t = h^2$}\label{table:3}
\end{table}
\end{center}
\begin{center}
\begin{table}[!h]
 \begin{tabular}{ |c||c|c|c|c|c|c|c|c| }
 \hline
J & $\mathcal{E}_1\times10$ & ${\rm eoc}_1$ & $\mathcal{E}_2$ & ${\rm eoc}_2$ & $\mathcal{E}_3$ 
& ${\rm eoc}_3$ & $\mathcal{E}_4\times10$ & ${\rm eoc}_4$ \\ 
 \hline
 \hline
 30 & 1.0152550 & - & 13.267420 & - & 3.8386880 & - & 2.5325450 & -   \\ 
 60 & 0.1996360 & 2.35 & 2.4120110 & 2.46 & 0.7548792 & 2.35 & 0.6895750 & 1.88  \\ 
120 & 0.0433905 & 2.20 & 0.5122927 & 2.24 & 0.1665889 & 2.18 & 0.1772575 & 1.96  \\ 
240 & 0.0100380 & 2.11 & 0.1179202 & 2.12 & 0.0391495 & 2.09 & 0.0449715 & 1.98  \\ 
\hline
\end{tabular}
\caption{$\alpha = 0.1, \Delta t = 0.5h$}\label{table:4}
\end{table}
\end{center}

\begin{figure}[h]
\centering
\subfigure{\includegraphics[width = 0.32\textwidth]{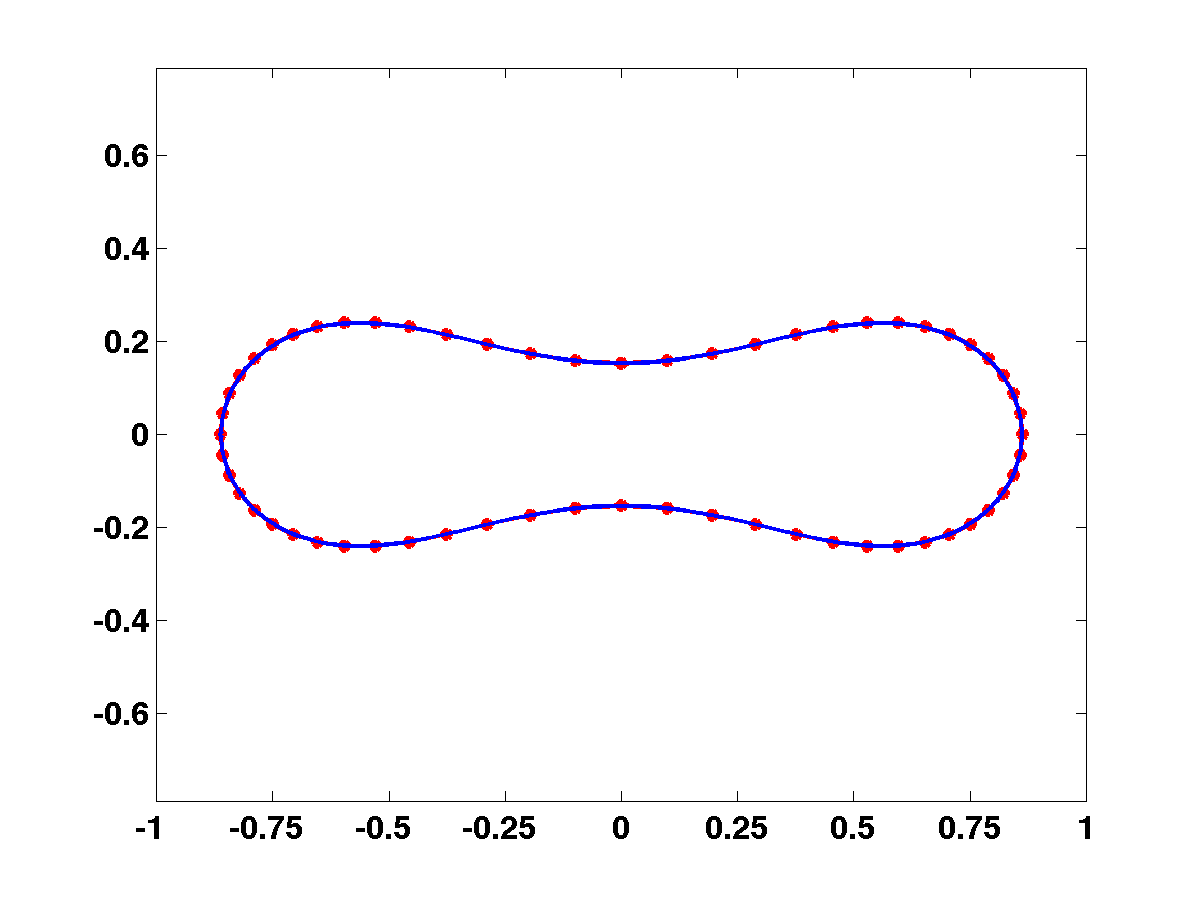}} 
\subfigure{\includegraphics[width = 0.32\textwidth]{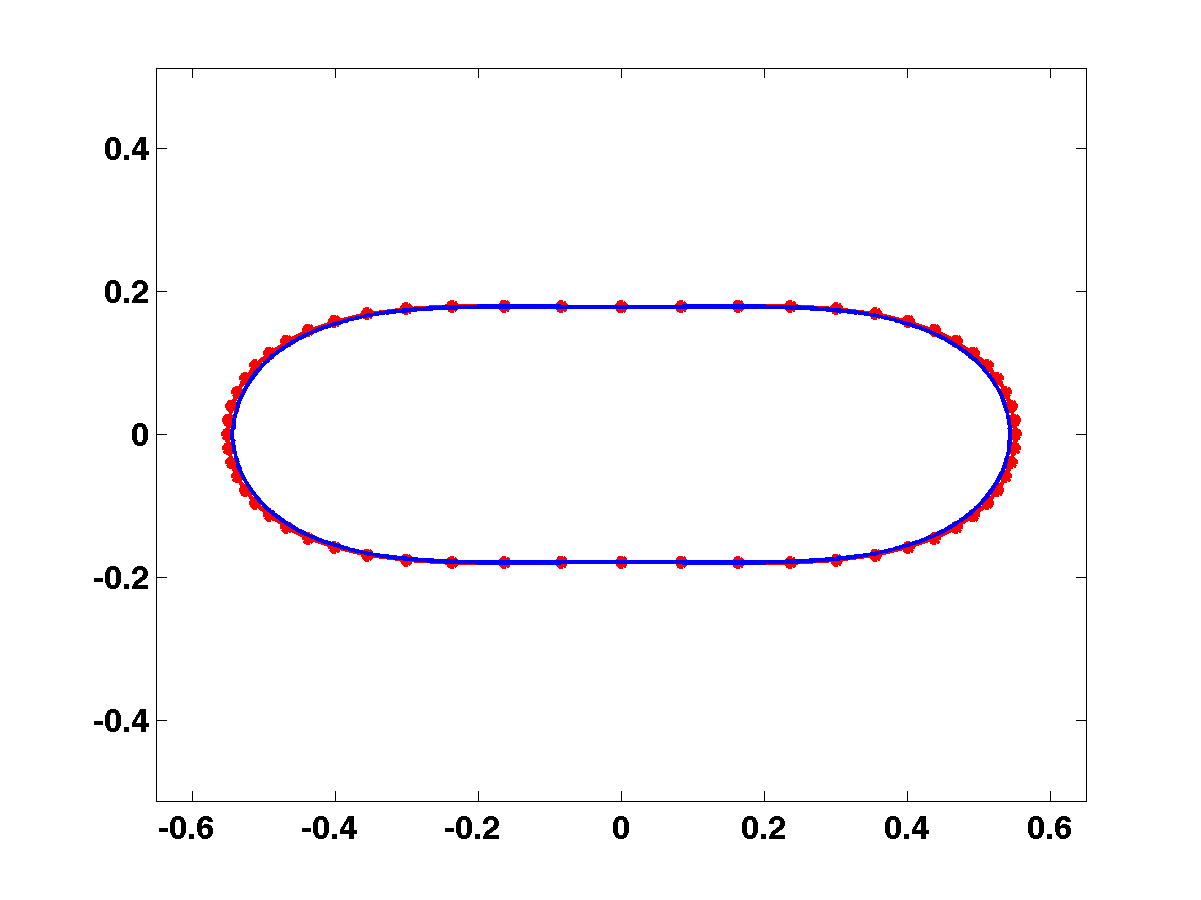}} 
\subfigure{\includegraphics[width = 0.32\textwidth]{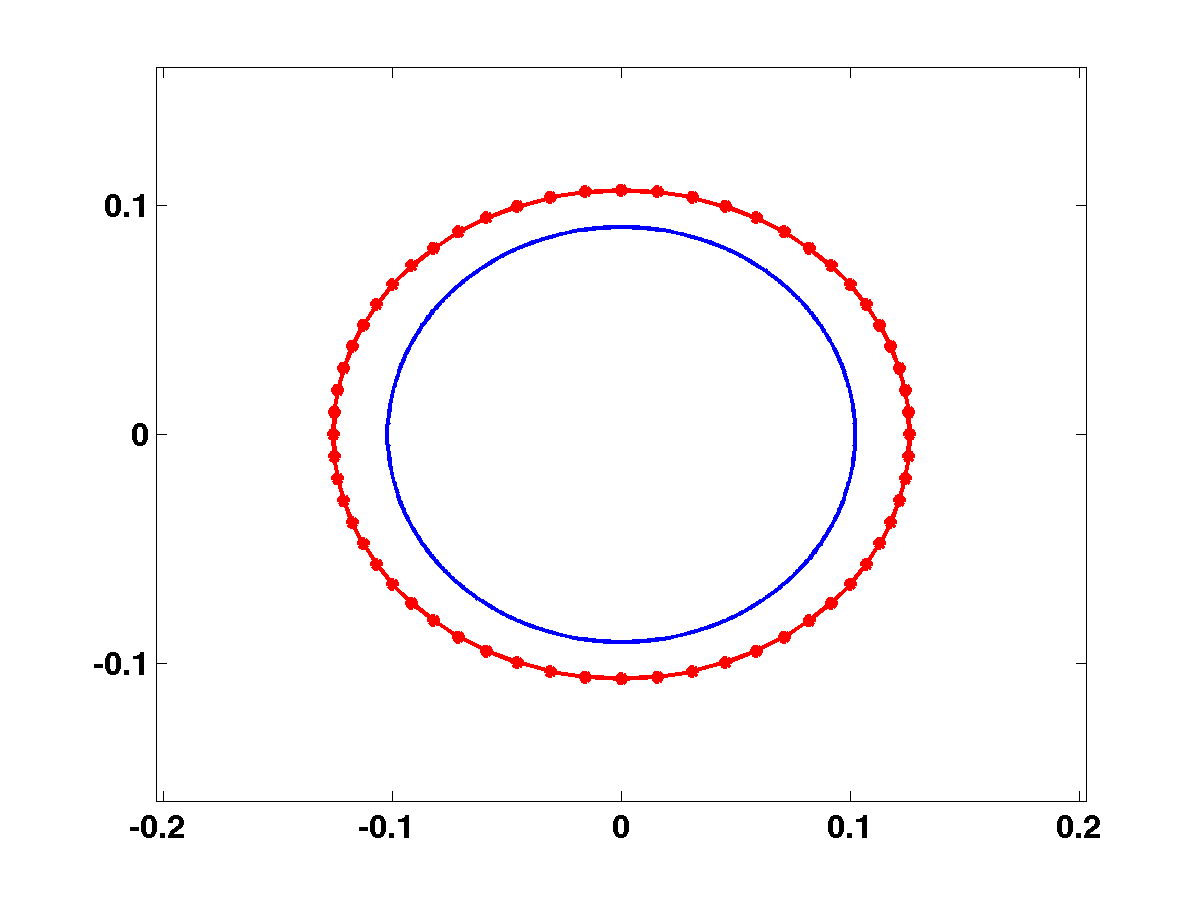}} \\
\subfigure{\includegraphics[width = 0.32\textwidth]{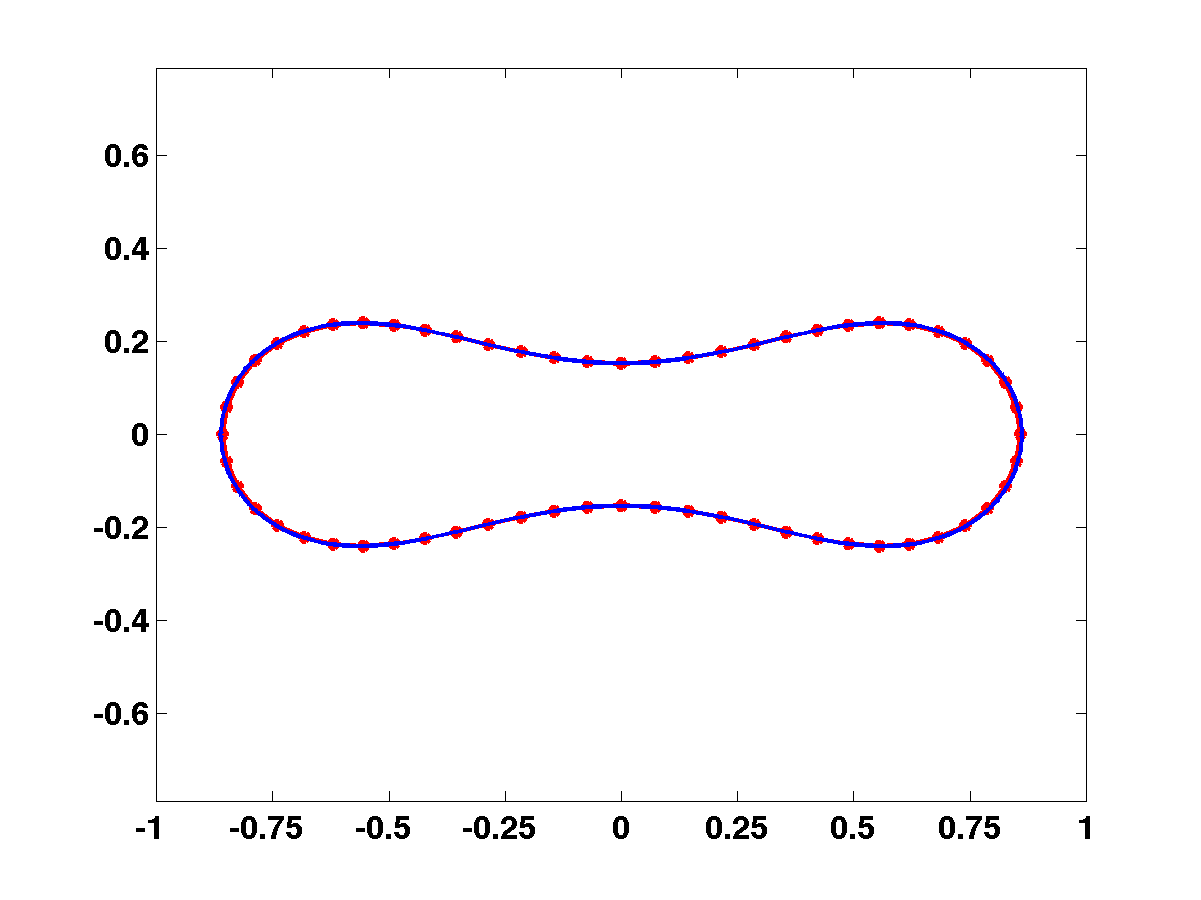}} 
\subfigure{\includegraphics[width = 0.32\textwidth]{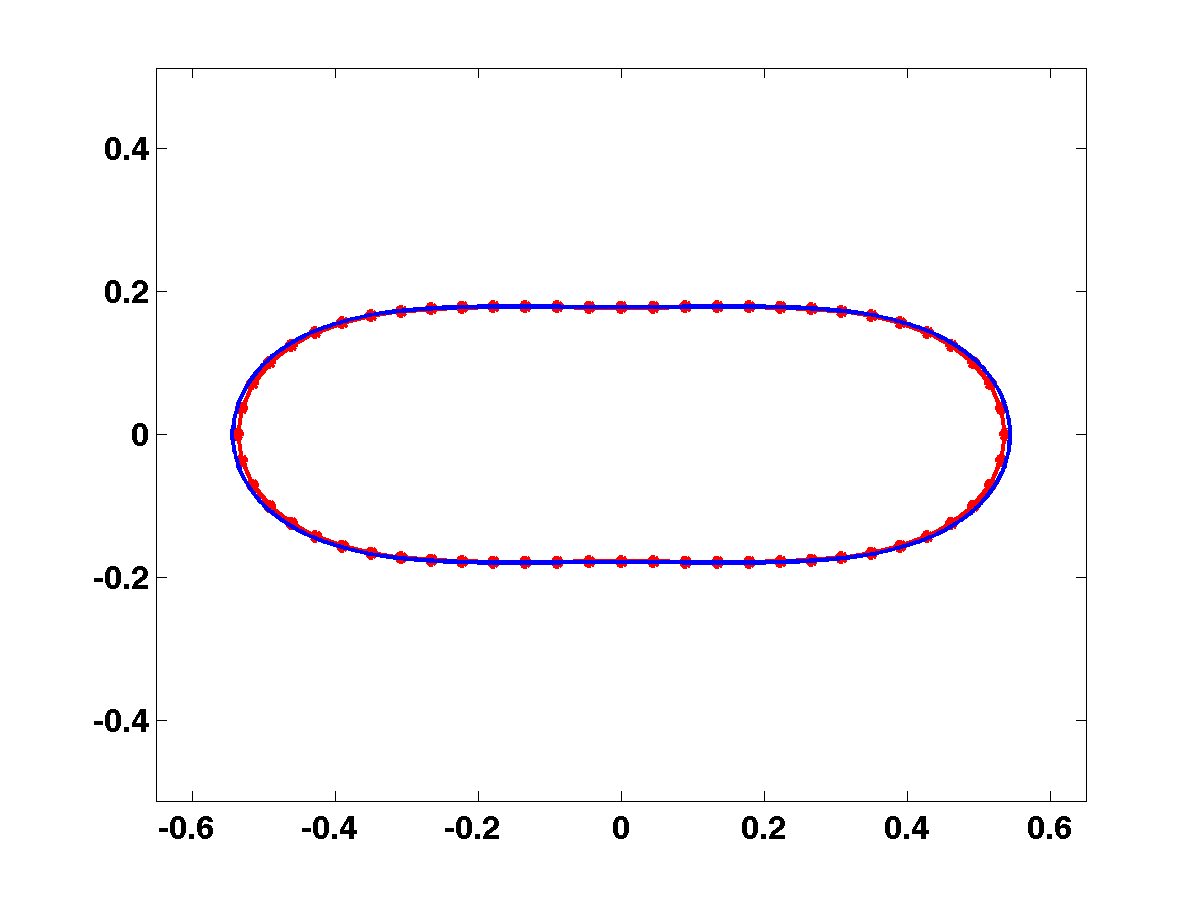}} 
\subfigure{\includegraphics[width = 0.32\textwidth]{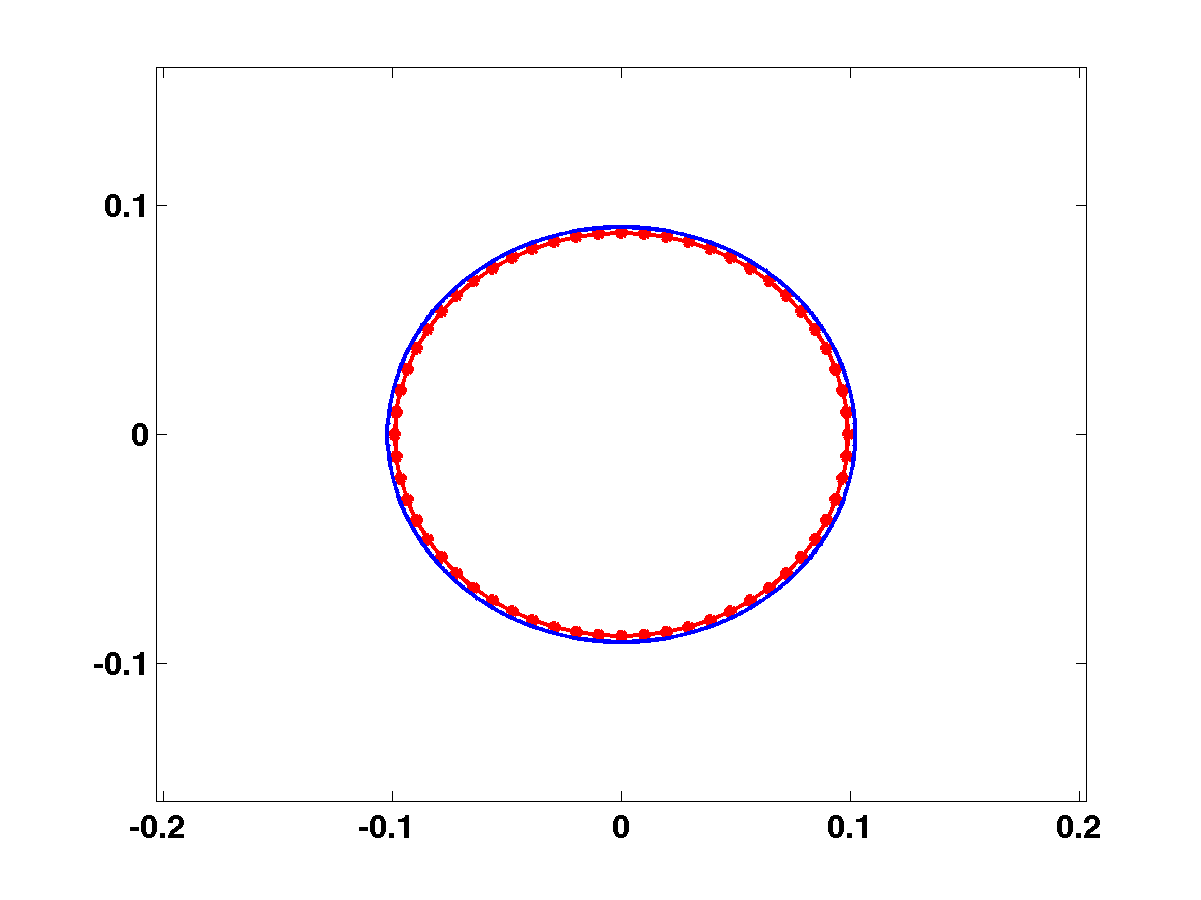}} 
\caption{$\vec{X}(\rho,t)$ at times $t=0.05,~0.1$ and $0.15$ with $\alpha=1$ 
(upper plots) and $\alpha =0.1$ (lower plots), $J=60$ red lines, $J=1200$ blue lines.}
\label{f:ex2_1}
\end{figure}

\begin{figure}[h]
\centering
\subfigure{\includegraphics[width = 0.32\textwidth]{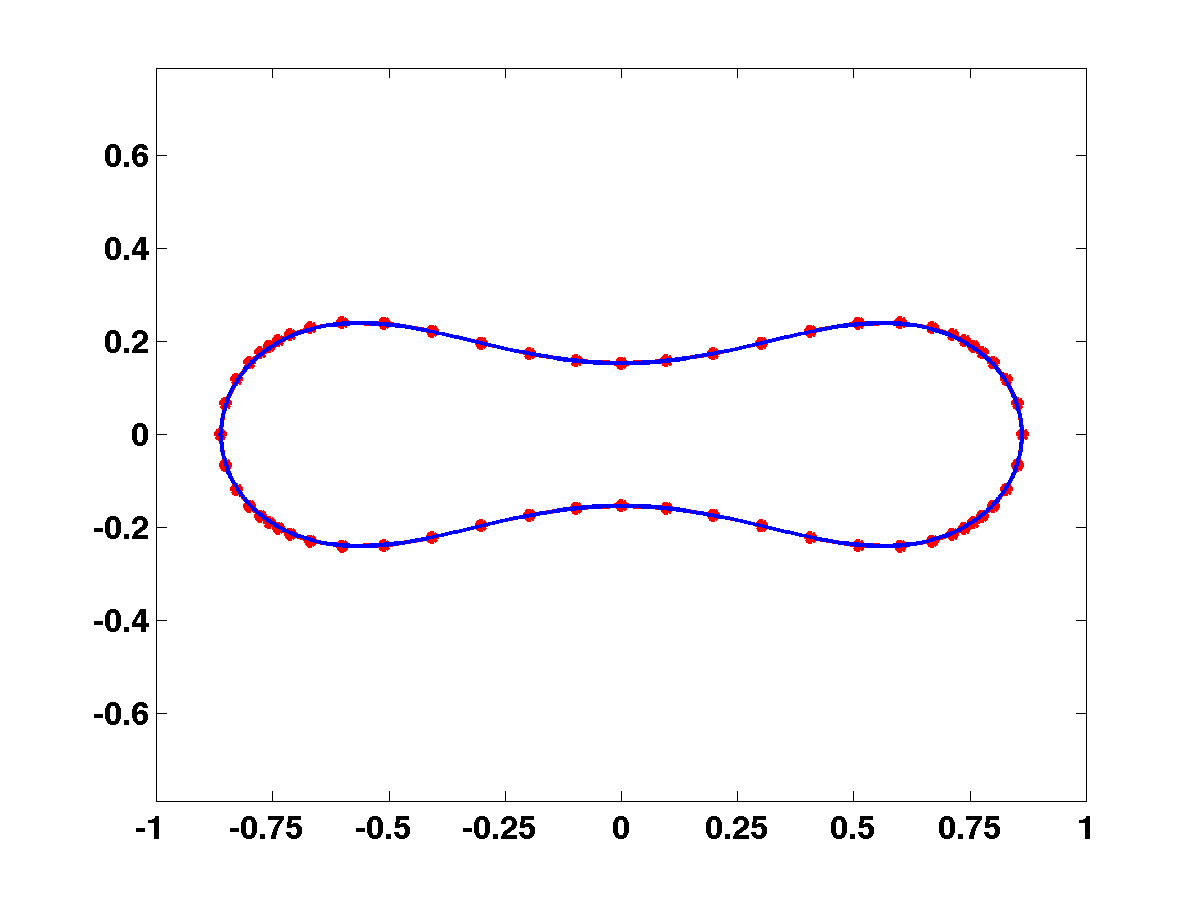}} 
\subfigure{\includegraphics[width = 0.32\textwidth]{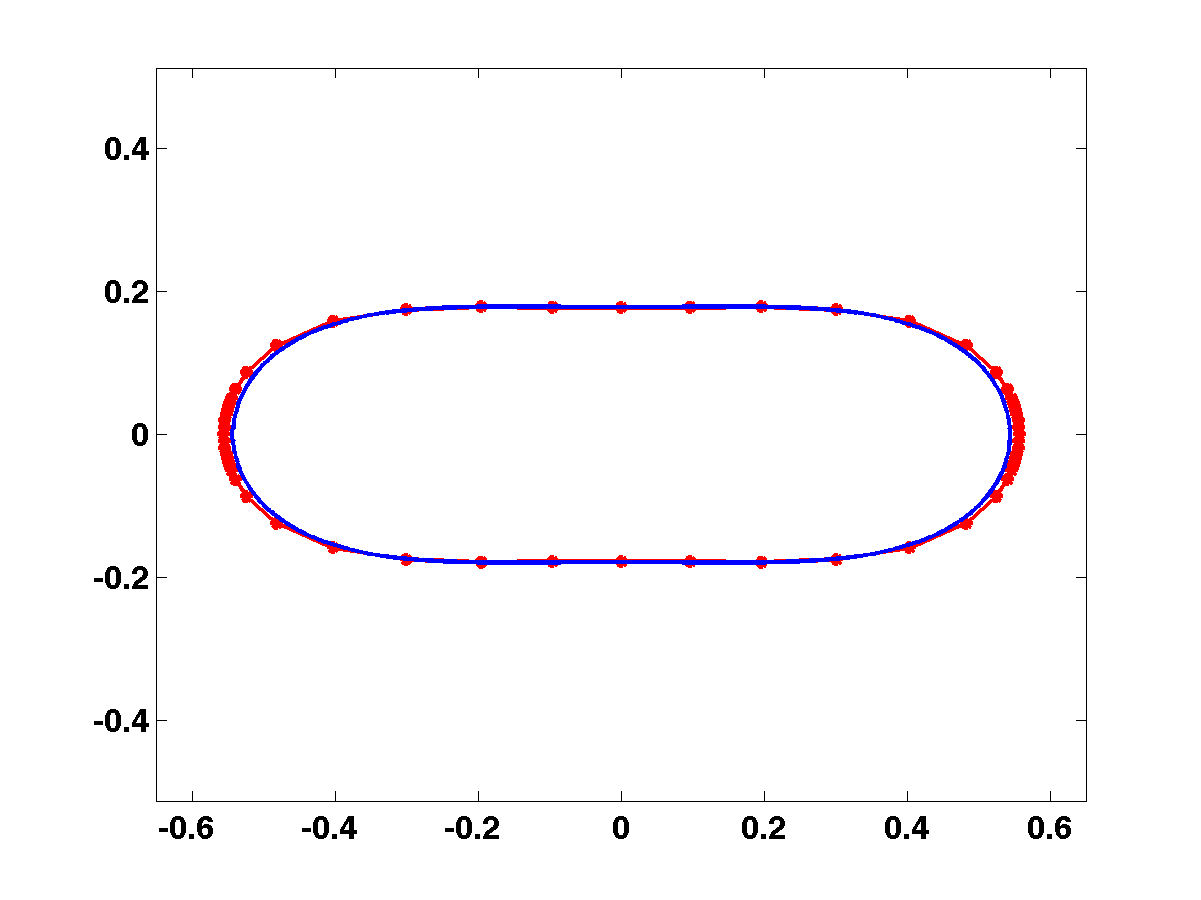}} 
\subfigure{\includegraphics[width = 0.32\textwidth]{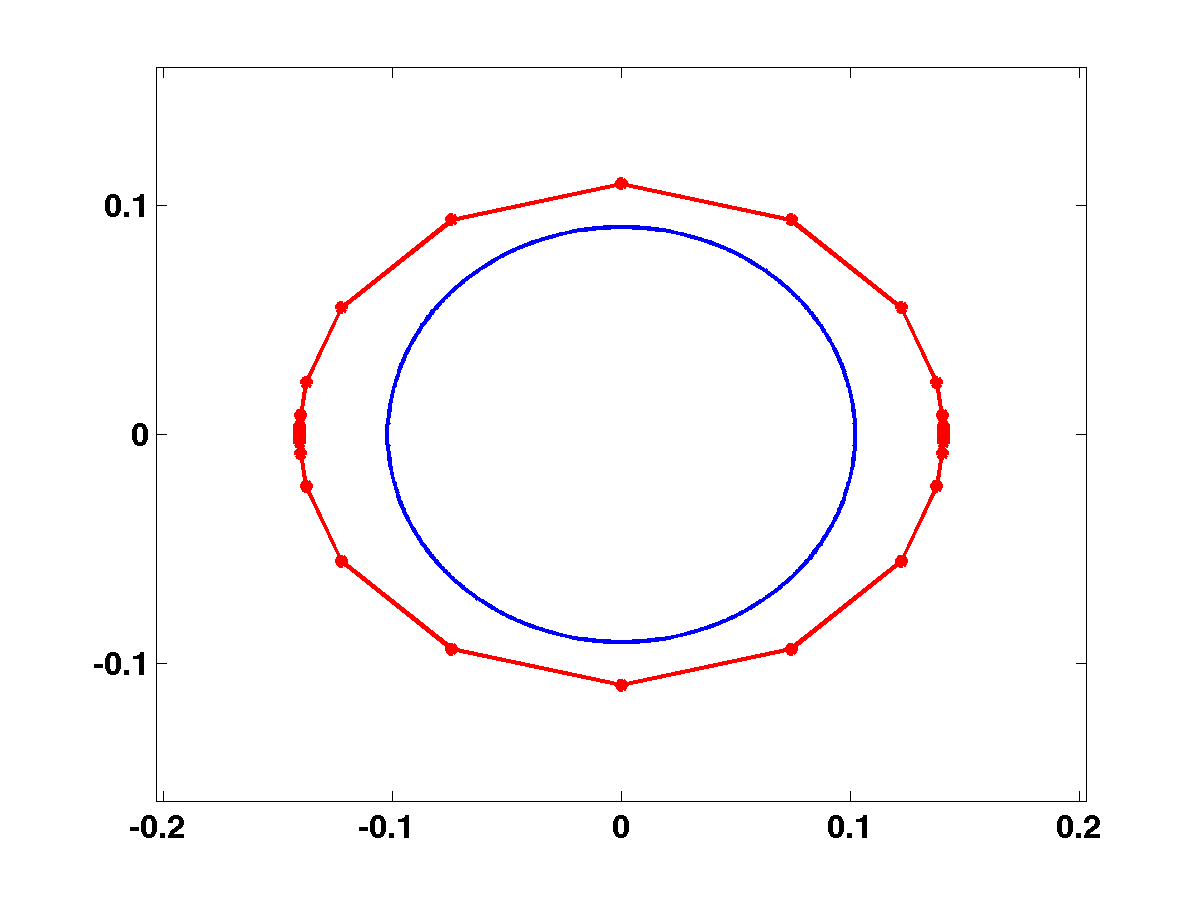}} 
\caption{$\vec{X}(\rho,t)$ at times $t=0.05,~0.1$ and $0.15$. The red line are computed using the scheme in \cite{PozziS15} with $J=60$ and the blue lines are computed using  (\ref{G1wbn}), (\ref{G2wbn}) with $J=1200$ and $\alpha=1$.}
\label{f:ex2_ps}
\end{figure}

In our second example we present numerical results that highlight 
the advantage of using $\alpha\ll 1$ for practical choices of $J$.
We set $f(w)=0.5w^2$, $g\left(v,w\right)=0$,
$$
\vec{x}(\rho,0)=\left(\begin{array}{cl}\cos(2\pi \rho)\\
(0.9\cos^2(2\pi \rho)+0.1)\sin(2\pi \rho) \end{array}\right), \quad 
w(\rho,0)=\sin(6\pi \rho), \quad \rho \in \Ip.
$$ 
For $J=1200$ the results obtained by plotting $\vX(\rho,t)$, for $t\in [0,0.15]$, 
with $\alpha=1$ and $\alpha=0.1$ are visually indistinguishable. 
However when $J=60$ at $t=0.15$ there is a marked difference between the two solutions. 
This can be seen in  
Figure \ref{f:ex2_1} where we display $\vX(\rho,t)$  (red lines) 
at times $t=0.05,~0.1$ and $0.15$ for $\alpha =1$ (upper plots) and $\alpha=0.1$ (lower plots). 
In these results we include $\vX(\rho,t)$ (blue lines), obtained by setting $J=1200$ and 
$\alpha=1$, to act as a comparison to the `true solution'. We see that the additional 
tangential motion of the nodes obtained by setting $\alpha=0.1$ results in a better 
approximation at $t=0.15$ to the `true solution' %with $J=1200$ 
than the approximation obtained by setting $\alpha=1$.  
The importance of the tangential motion of the nodes can be further seen in Figure \ref{f:ex2_ps} where we display 
$\vX(\rho,t)$ (red lines) at times $t=0.05,~0.1$ and $0.15$ for the scheme in \cite{PozziS15} with $J=60$,
in which there is no tangential motion. Again 
we include %$\vX(\rho,t)$ (blue lines), obtained by setting $J=1200$ and 
%$\alpha=1$ in (\ref{G1wbn}), (\ref{G2wbn}) as an approximation to the 
`true solution' (blue lines) as in Figure \ref{f:ex2_1}.
The absence of tangential motion leads to a severe accumulation of the nodes at $t=0.15$, however 
as the value of $J$ increases this accumulation effect becomes less pronounced.

\begin{figure}[h]
\centering
\subfigure{\includegraphics[width = 0.32\textwidth]{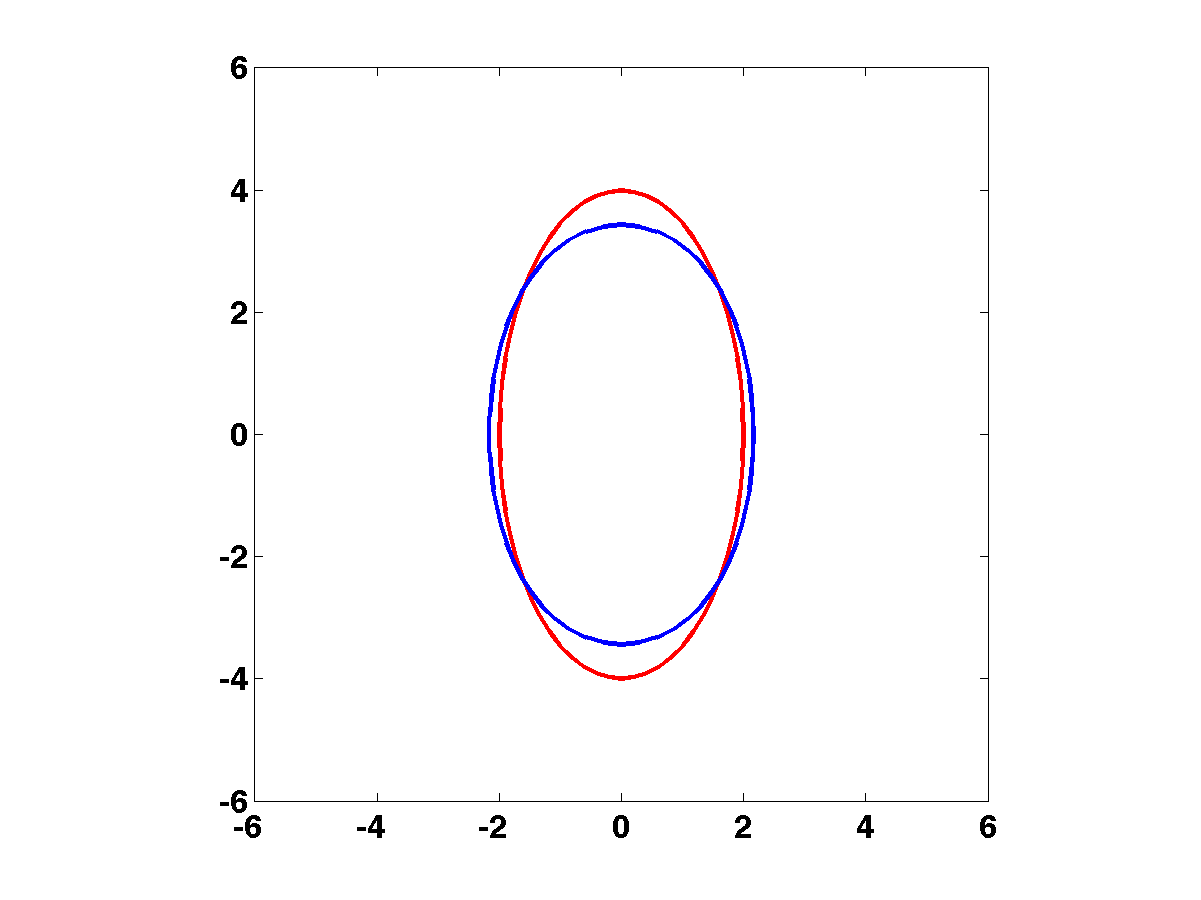}} 
\subfigure{\includegraphics[width = 0.32\textwidth]{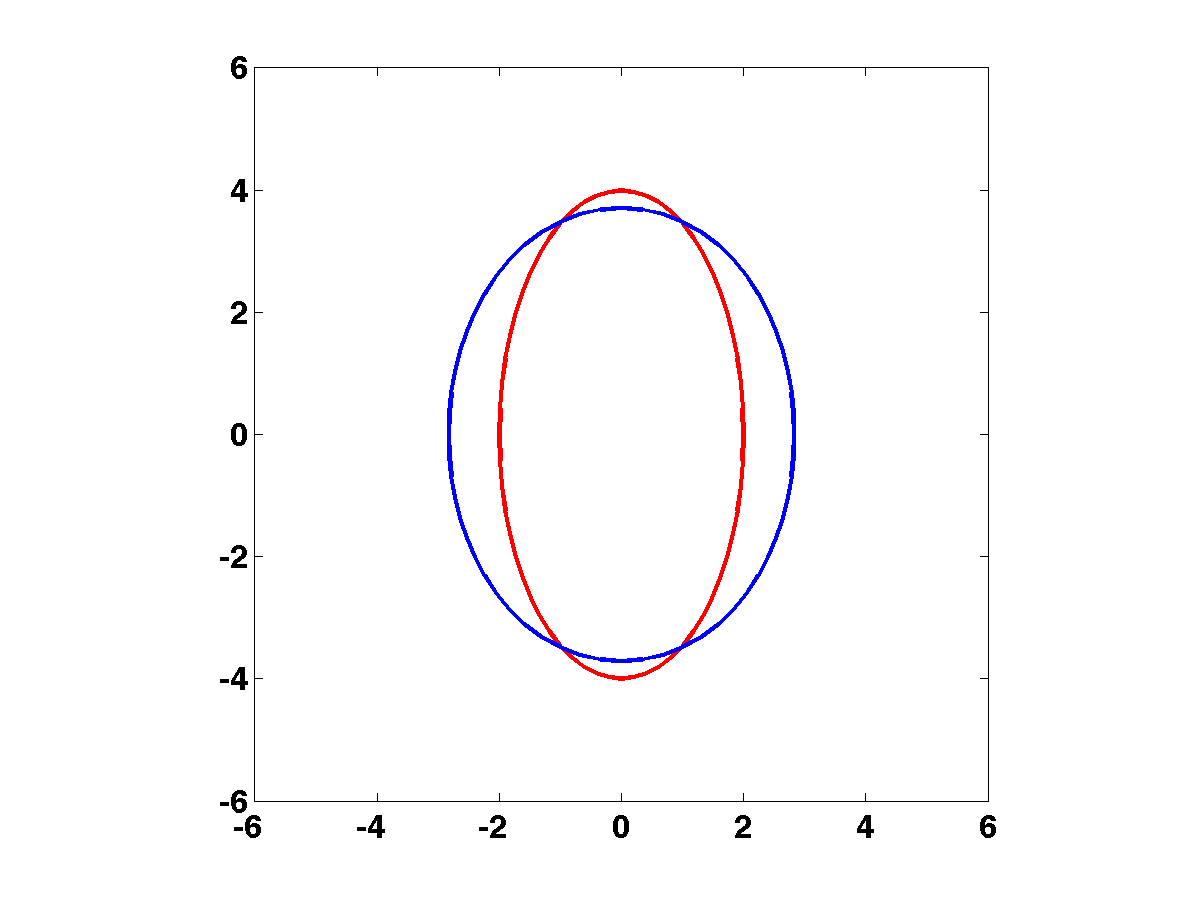}}  
\subfigure{\includegraphics[width = 0.32\textwidth]{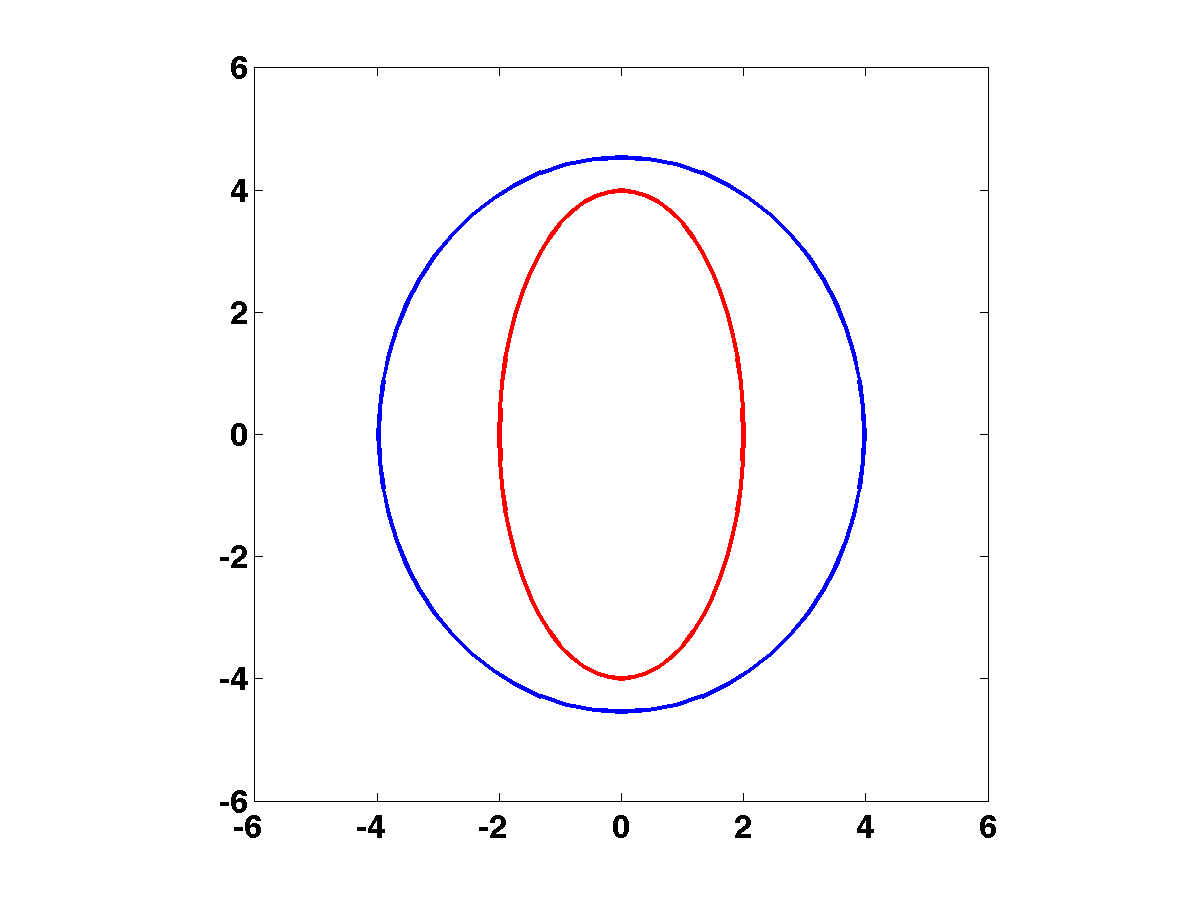}} 
\caption{The grain boundary $\vX(\rho,t)$  at times $t=1.5,~3.5,~7.0$ (blue lines) and at $t=0$ (red lines).}
\label{f:digm1}
\end{figure}

\begin{figure}[h]
\centering
\subfigure{\includegraphics[width = 0.28\textwidth]{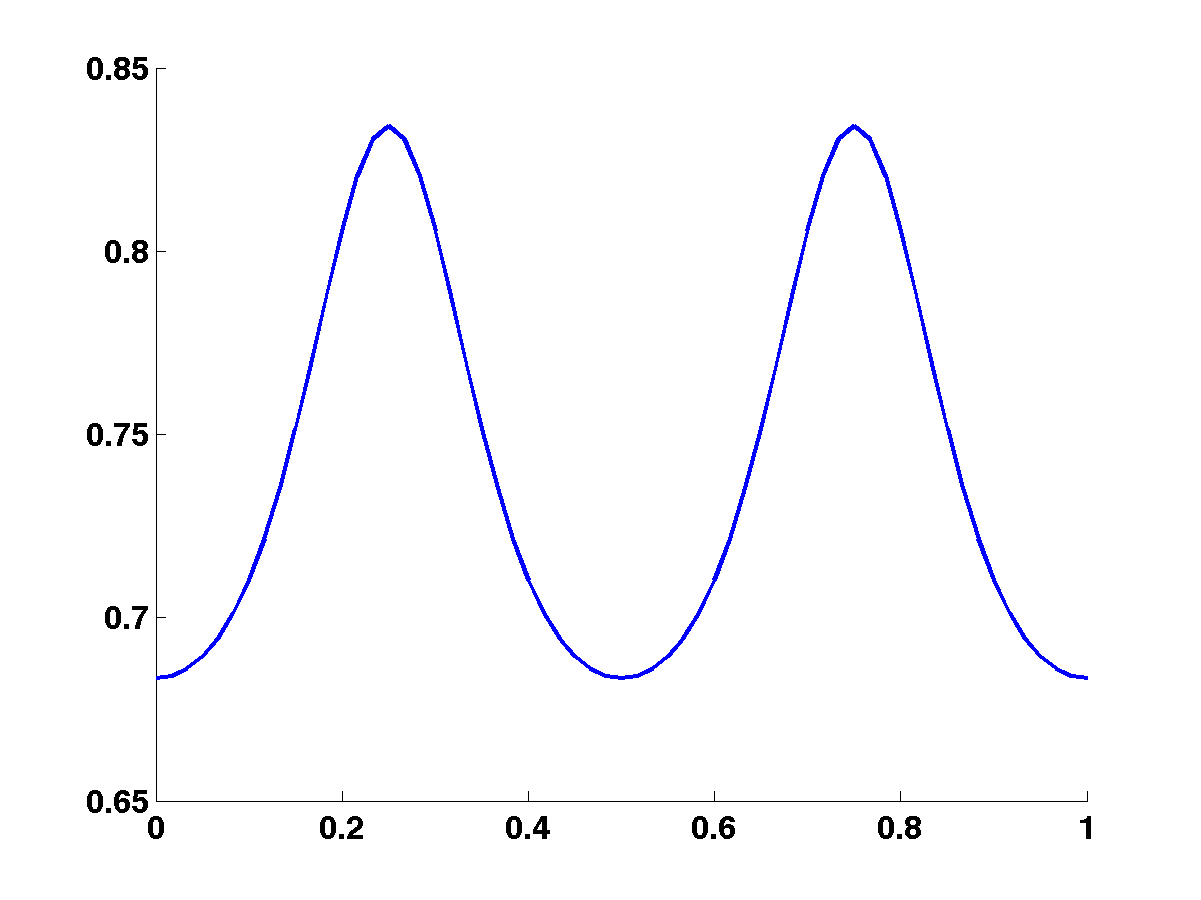}} 
\subfigure{\includegraphics[width = 0.28\textwidth]{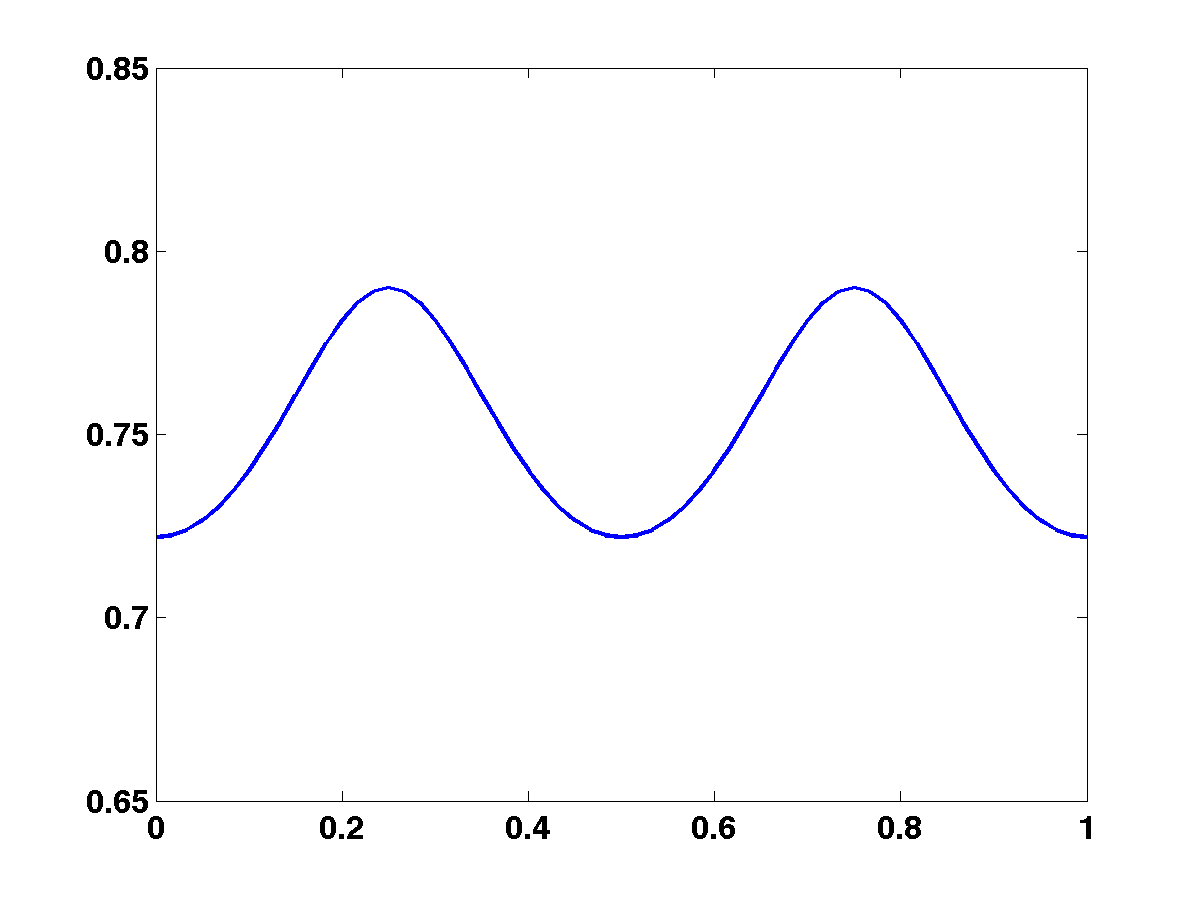}}  
\subfigure{\includegraphics[width = 0.28\textwidth]{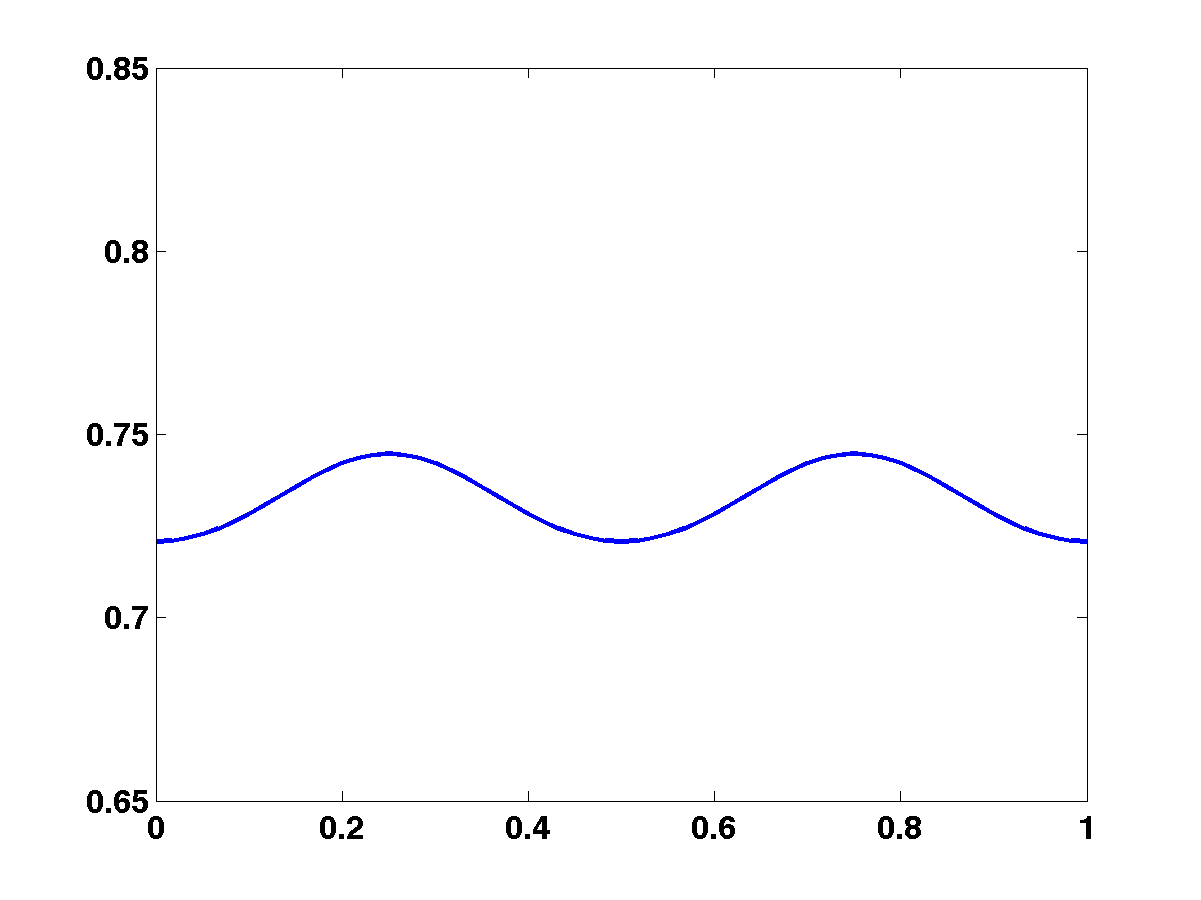}} 
\caption{The concentration $W(\rho,t)$ at times at times $t=1.5,~3.5,~7.0$.}
\label{f:digm2}
\end{figure}

We conclude our numerical results with a simulation of diffusion induced grain boundary motion. 
%In \cite{MayerS99} it is stated that 
With $\vnu$ being the inward normal, 
physically meaningful choices for 
$f(w)$ and $g(v,w)$ in (\ref{Gamma1},b) are $f(w)=-w^2$ and 
$g(v,w)=v\,w-(C-w)$,
where  $C$ is the concentration of solute in the vapour. Note that these choices satisfy
(\ref{fgreg},b).
We consider the experimental setup in which the parametrisation of the initial grain boundary 
and the initial concentration of solute are given by
$$
\vec{x}(\rho,0)=\left(\begin{array}{cl}2\cos(2\pi \rho)\\
4\sin(2\pi \rho) \end{array}\right), \quad 
w(\rho,0)=0, \quad \rho \in \Ip
$$ 
respectively, while  $C=1$. In the simulations presented we set $\alpha =0.1$ and $J=60$. 
Figure \ref{f:digm1} displays the grain boundary, $\vX(\rho,t)$, 
(blue line) at times $t=1.5$, $t=3.5$ and $t=7.0$ together with its initial position (red lines). 
At $t=1.5$ we see an example of bi-direction motion occurring, 
see \cite{ElliottS03}, \cite{GarckeS04}. 
In Figure \ref{f:digm2} we plot the concentration $w(\rho,t)$ at $t=1.5$, $t=3.5$ and $t=7.0$.

\section*{Acknowledgements}
VS gratefully acknowledges the supported by the Engineering and Physical Sciences Research Council, UK grant (EP/J016780/1) and the Leverhulme Trust  Research Project Grant (RPG-2014-149). 

\bibliographystyle{siamplain}

%\bibliography{refsjwb}

%\end{document}

\end{document}